# INDUCED GELATION IN A TWO-SITE SPATIAL COAGULATION MODEL[1]

By Rainer Siegmund-Schultze and Wolfgang Wagner

*Technische Universität Berlin and Weierstrass Institute for Applied Analysis and Stochastics*

A two-site spatial coagulation model is considered. Particles of masses $m$ and $n$ at the same site form a new particle of mass $m+n$ at rate $mn$. Independently, particles jump to the other site at a constant rate. The limit (for increasing particle numbers) of this model is expected to be nondeterministic after the gelation time, namely, one or two giant particles randomly jump between the two sites. Moreover, a new effect of induced gelation is observed—the gelation happening at the site with the larger initial number of monomers immediately induces gelation at the other site. Induced gelation is shown to be of logarithmic order. The limiting behavior of the model is derived rigorously up to the gelation time, while the expected post-gelation behavior is illustrated by a numerical simulation.

**1. Introduction.** We consider a minimal spatial model of coagulation with the gelling kernel $mn$ describing the coagulation rate and an independent motion of particles between two sites 0 and 1 which happens at a constant rate for each particle irrespective of its mass. So the basic model is a continuous-time Markov jump process

$$\{\xi_{t,m,j}\}, \qquad t\geq 0,\ m=1,2,\ldots,\ j=0,1.$$

A state of the process is given by two sequences of nonnegative integers, both with only a finite number of nonzero entries, describing the numbers of particles of different masses at the two sites. As for the initial distribution, we consider the deterministic state given by the two length-1 sequences $\{N_0\}, \{N_1\}$ meaning a start with $N_0, N_1$ monomers at the two sites. Denote

Received May 2004; revised August 2005.
[1]Supported by Deutsche Forschungsgemeinschaft Grant BA 922/7-1.

the total number of monomers by

$$N = N_0 + N_1.$$

Given a state $(\xi_{m,j})$, the jump rate (depending on the scaling parameter $N$) has the form

$$\tfrac{1}{2}N(N-1) + \kappa N \sum_{m,j} \xi_{m,j}.$$

The transition is defined as follows: Either, with probability

$$\kappa \sum_{m,j} \xi_{m,j} \left( \tfrac{1}{2}(N-1) + \kappa \sum_{m,j} \xi_{m,j} \right)^{-1},$$

a randomly selected particle (an aggregation of monomers) jumps to the opposite site, or with the remaining probability two monomers are randomly chosen and the corresponding particles coalesce provided they are not the same particle and they are at the same site. This means that an unordered pair of different particles of masses $m$ and $n$ at the same site coalesces at rate $mn$. Hence our model is a (minimal) spatial generalization of the Marcus–Lushnikov stochastic coalescent $\mathrm{ML}^{(N)}(t)$ with monomeric initial condition [12, 13].

It is clear that we can describe the evolution as well as a colored graph evolution (with two colors): Start with $N_0$ isolated red nodes and $N_1$ isolated green nodes. Add at constant rate 1 an edge for each unordered pair of nodes, where edges between differently colored nodes are instantly removed as well as multiple edges, and let (maximal) connected components change their color at constant rate $\kappa N$.

We investigate the limiting behavior of this model for $N$ tending to infinity while

$$(1) \qquad N_1/N_0 \to \lambda,$$

where $\lambda$ is a given parameter. The expected behavior is notably more complex compared to the well-known case where all particles live at the same site (*mean-field* model or homogeneous model). The mean-field model has a deterministic limit, described by a modified Smoluchowski system of differential equations (with the modification being urged by the effect of *gelation*, meaning that, after *rescaling time* by $t \to t/N$, at gelation time $t_{\mathrm{gel}} = 1$ a giant particle of mass comparable to the total mass suddenly appears in a phase transition and modifies the evolution). The more general model considered here is expected to have a stochastic limit after gelation time, caused by the fact that the motion of (at most two) giant particles is no longer governed by the law of large numbers. Moreover, a new effect of *induced gelation* is observed: The gelation happening at the site with the larger initial concentration of monomers immediately induces a gelation at the other site. We give a rigorous derivation of the limiting behavior up to the gelation time which is based mainly on the cutting techniques developed in [14, 15].

We emphasize that a mathematically rigorous derivation of the limiting behavior after gelation in our model can be expected to be quite difficult since the existing techniques either are based on a martingale approach and rely heavily upon the finiteness of the expected particle size per monomer [14], which breaks down at the gelation time, or they make use of random graph techniques [2] which admit no (or at least no immediate) generalization to the spatial model. Also, as will be seen later, the speed of convergence to the limit process is very slow, namely of logarithmic order in the total particle number, because induced gelation is an effect of logarithmic effectiveness. Hence, as far as *post-gelation* behavior is concerned (Section 4), this account is of mainly phenomenological nature and will deliver evidence for the expected phenomena by numerical simulation. In the rigorous part we use the notions "gelation time" and "gelation" in a weaker sense, meaning the explosion of the second moment.

**2. Limiting behavior of the spatial model before gelation.** Here we study the two-site model at pre-gelation times. The parameter $\lambda$ indicates the dependence on the initial condition (1).

THEOREM 1. *Let $\kappa > 0$. Then there is a nonrandom time $t_{\mathrm{gel}}^\lambda > 0$, such that the rescaled quantities $N^{-1}\xi_{t/N,m,j}$ of our model with asymptotic parameter $\lambda$ do uniformly converge in probability to deterministic quantities $c_{t,m,j}^\lambda$ fulfilling a two-site analogue of Smoluchowski's coagulation equation:*

$$\dot{c}_{t,m,j}^\lambda = \tfrac{1}{2} \sum_{0<n<m} n(m-n) c_{t,n,j}^\lambda c_{t,m-n,j}^\lambda - \sum_{n>0} mn c_{t,m,j}^\lambda c_{t,n,j}^\lambda$$
(2)
$$- \kappa c_{t,m,j}^\lambda + \kappa c_{t,m,1-j}^\lambda$$

*with initial condition*

$$c_{0,1,0}^\lambda = \frac{1}{1+\lambda}, \qquad c_{0,1,1}^\lambda = \frac{\lambda}{1+\lambda}, \qquad c_{0,m,j}^\lambda = 0, \qquad m > 1.$$

*Moreover, $t_{\mathrm{gel}}^\lambda$ is characterized by the fact that the two quantities*

$$\sigma_{t,j}^\lambda = \sum_{m>0} m^2 c_{t,m,j}^\lambda, \qquad j = 0, 1,$$
(3)

*reach infinity for the first time.*

By a *strong maximal positive solution* of the system of ODE (2) we understand a solution defined in an interval $[0, \tilde{t})$ which obeys the properties:

(a) $0 \leq c^\lambda_{t,m,j}$ for all $m \in \mathbb{N}, j = 0, 1, t \in [0, \tilde{t})$,
(b) the second moments (3) exist for $t \in [0, \tilde{t})$ and fulfill

$$\int_0^T \sigma^\lambda_{s,j}\, ds < \infty, \qquad T < \tilde{t},\ j = 0, 1,$$

(c) $0 < \tilde{t} \leq +\infty$ is the largest value with properties (a), (b).

THEOREM 2. *The two-site Smoluchowski equation* (2) *has a unique strong maximal positive solution with* $\tilde{t} = t^\lambda_{\mathrm{gel}}$.

Next, let us consider the differential equation describing (as will be seen later, at the end of the proof of Theorem 2) the behavior of the second moments at both sites at pre-gelation times

(4) $$\dot\sigma^\lambda_{t,j} = (\sigma^\lambda_{t,j})^2 + \kappa(\sigma^\lambda_{t,1-j} - \sigma^\lambda_{t,j}),$$

with initial condition

$$\sigma^\lambda_{0,0} = (1+\lambda)^{-1}, \qquad \sigma^\lambda_{0,1} = \lambda(1+\lambda)^{-1}.$$

Equation (4) is autonomous and extends in an obvious way the well-known equation $\dot\sigma = \sigma^2$ valid in the pre-gelation regime for the mean-field model. Writing $x(t)$ for $\sigma^\lambda_{t,0}$ and $y(t)$ for $\sigma^\lambda_{t,1}$, we get the following system of two ODEs:

(5) $$\dot x = x^2 + \kappa(y-x), \qquad \dot y = y^2 + \kappa(x-y)$$

with

(6) $$x_0 = (1+\lambda)^{-1} \geq 0, \qquad y_0 = \lambda(1+\lambda)^{-1} \geq 0.$$

THEOREM 3. *If* $\kappa > 0$, *then there is a positive* $t_{\mathrm{gel}} \leq 2$ *such that* $x, y$ *are smooth in* $[0, t_{\mathrm{gel}})$ *and we have* $x, y \to +\infty$ *as* $t \to t_{\mathrm{gel}}$ *simultaneously.*

THEOREM 4. *If* $\kappa > 0$ *and* $x_0 > y_0$, *then* $x > y$ *for all* $t \in [0, t_{\mathrm{gel}})$ *and* $x \sim (t_{\mathrm{gel}} - t)^{-1}$ *as* $t \to t_{\mathrm{gel}}$ *while* $y \sim -\kappa \log(t_{\mathrm{gel}} - t)$.

Hence (except for $\kappa = 0$) *both* second moments (expected particle masses at site $j$ per monomer) $\sigma^\lambda_{t,j}, j = 0, 1$, become infinite simultaneously (this is trivial for $\lambda = 1$), an effect we call *induced gelation* in view of the mechanism lying behind this: If $\lambda \neq 1$, the immigration shortly before gelation time of large particles from the site with the larger initial number of particles (i.e., site 0 in the case $0 \leq \lambda < 1$) forces a gelation on the other site. The numerical studies with particle numbers up to $10^9$ show that this induced gelation exhibits some special features which distinguish it clearly from the "genuine" gelation taking place at the first site. In fact, for $\lambda < 1$ the growth of $\sigma_{t,0}^\lambda$ is like $(t - t_{\text{gel}})^{-1}$ while the growth of $\sigma_{t,1}^\lambda$ is only of logarithmic order $\kappa \log(t - t_{\text{gel}})^{-1}$. Moreover, the proof of Theorem 4 indicates that, approaching $t_{\text{gel}}$, the growth of $\sigma_{t,1}^\lambda$ is finally completely caused by particles with large masses immigrating from site 0, while the "genuine" growth by coagulation at site 1 can be asymptotically neglected.

On the other hand, it is clear from the logarithmic character of the induced gelation that it can be expected to be practically perceptible only for particle numbers even the logarithm of which is large enough. Take, for instance, $10^9$ particles at site 0 and 0 particles at site 1 with $\kappa = 0.3$. As is known from the one-site model, at gelation time the size of the giant particle can be expected to be about $10^{(2/3)9} = 10^6$. So in fact the mass of immigrating particles up to gelation time is bounded by this quantity. Therefore the maximum increase of the (empirical) value for $\sigma_{t,1}^\lambda$ is most probably less than $6\kappa \log 10 \approx 4$. That is an estimate of the maximal increase induced by site 0 shortly before gelation happens there. The gelation at site 1 has to be accomplished now by coagulation of particles at site 1, since the immigration of massive particles from site 0 rather rapidly stops after gelation, because these are rapidly consumed by the giant particle. The estimated remaining time it takes for site 1 to reach gelation is now about $1/4$, meaning that in fact even for $10^9$ particles a perceptible delay occurs until the induced gelation actually takes place. Making $\kappa$ larger does not help too much, because then the tendency to reach equilibrium between both sites is stronger, such that it becomes harder to distinguish the induced from the "genuine" gelation, which would have taken place anyway at site 1, if migration was stopped shortly before $t_{\text{gel}}$. Nonetheless, as we will see (interpreting the numerical studies), it is for the given choice of parameters not difficult to recognize the induced gelation by some characteristic properties appearing after $t_{\text{gel}}$.

**3. Properties of the mean-field model.** Consider the case $\kappa = 0, N = N_0$, that is, the mean-field model. Then the graph is completely colored red for all times and we are in the situation considered in great detail in the monograph [2]. There the random graph model $\mathcal{G}(N, p)$ is studied which describes a graph with $N$ nodes, where for each pair of nodes an edge is present with probability $p$, completely independent from the presence or nonpresence of the other edges. Since the process described above connects an unordered pair of different monomers at constant rate 1, the probability of the corresponding edge to be present at time $t$ is $1 - e^{-t}$. Hence the distribution of $\{\xi_{t,m,0}\}_{m>0}$ is given by the joint distribution of maximally connected components of sizes $m = 1, 2, 3 \ldots$ in the random graph $\mathcal{G}(N, 1 - e^{-t})$.

As is well known in both the coagulation interpretation (cf. the comprehensive paper [1]) and the random graph setup [2, 5], the time-rescaled ($t \to t/N$) model undergoes at $t = 1$ (i.e., at time $1/N$ in the old scale) a phase transition called gelation, which is characterized by the occurrence of a single "giant particle" of mass $\xi_t^{\max} = \Theta(N)$ [comparable to $N$ in the sense that $(\xi_t^{\max})^{-1} = O(N^{-1})$] for $t > 1$. Up to this gelation time $t_{\text{gel}} = 1$ the rescaled particle numbers $N^{-1}\xi_{t,m,0}$, $m = 1, 2, \ldots$, tend in probability to deterministic quantities $c_{t,m}$ fulfilling the Smoluchowski system of differential equations

$$\dot{c}_{t,m} = \tfrac{1}{2} \sum_{0<n<m} n(m-n) c_{t,n} c_{t,m-n} - m c_{t,m} \sum_{n>0} n c_{t,n} \tag{7}$$

which conserves the total (rescaled) mass 1, that is, $\sum_{m>0} m c_{t,m} \equiv 1$, $t \leq 1$. At later times the limiting behavior is *no longer* described by this set of equations. As expected, the Smoluchowski system no longer conserves the mass for $t > t_{\text{gel}}$, but it does not take into account the strong effect of the giant particle on the rest of the system. The actual limiting quantities of the stochastic coalescent differ from the solutions of (7), though sharing the property that total mass gets lost.

The proper system of equations is

$$\dot{c}_{t,m} = \tfrac{1}{2} \sum_{0<n<m} n(m-n) c_{t,n} c_{t,m-n} - m c_{t,m}. \tag{8}$$

Observe that in contrast to the infinite system of equations (7) the system (8) has the nice property that the finite subsystems obtained for $m \leq M$ are autonomous. It is obtained from (7) by replacing the quantity $\sum_{n>0} n c_{t,n}$ by 1, which simply means that particles are aggregated to others (including the giant one) at constant rate 1, representing the invariant rescaled total mass.

The solution to (8) can be given in explicit terms for the case $c_{0,1} = 1$, that is, the monodisperse initial condition:

$$c_{t,m} = t^{m-1} m^{m-2} e^{-tm} / m!. \tag{9}$$

The convergence of the rescaled stochastic coalescent to (9), that is, the solution of (8), was shown rigorously [2] for *all* $t \geq 0$.

The quantities $\{m c_{t,m}\}_{m>0}$, together with a mass of $\zeta_{t,\infty} := 1 - \sum_{m>0} m c_{t,m}$ assigned to the value $+\infty$, form the Borel distribution with parameter $t$, that is, the distribution of the total progeny of a Galton–Watson process with Poisson offspring of parameter $t$. Observe that for $t > 1$ the branching process is supercritical and hence the total offspring is infinite with positive probability, leading to a mass defect:

$$\zeta_{t,*} := \sum_{m>0} m c_{t,m} < 1,$$

reflecting gelation.

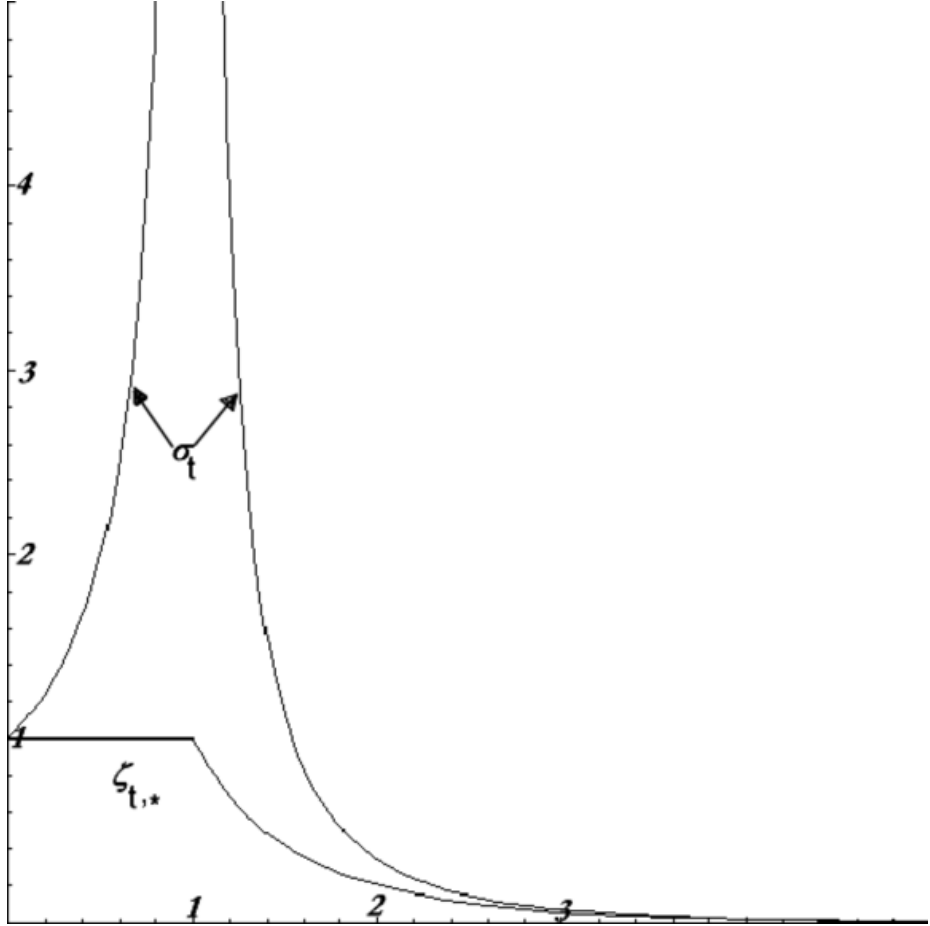

Fig. 1. *Plot of $\zeta_{t,*}$ and $\sigma_t$.*

While for the system (7) the second moment $\sigma_t := \sum_{m>0} m^2 c_{t,m}$ reaches infinity exactly at $t_{\text{gel}}$ and remains infinite, for the system (8) this second moment is infinite *only* at $t_{\text{gel}}$ (see Figure 1). In fact, $\sigma_t$ is just a reduced second moment, with the gel mass not being included.

The second moment has a more immediate interpretation as a *first moment*: It is the average mass of the particle where a randomly chosen monomer is aggregated into. So at gelation time this average mass becomes infinite, and it returns to finite values only because monomers "get lost" in the giant particle.

Let us look at the behavior of the system (8) in some more detail. For this sake it is convenient to change from the consideration of particle numbers to the *mass flow equation* which is obtained if we consider the quantities $\zeta_{t,m} := m c_{t,m}$ representing the mass concentrated in particles of size $m$ at time $t$. Now the evolution equation becomes

$$\dot{\zeta}_{t,m} = \tfrac{1}{2} m \sum_{0<n<m} \zeta_{t,n} \zeta_{t,m-n} - m \zeta_{t,m}. \tag{10}$$

We introduce the (*probability*) *generating function*

$$u(x,t) := \sum_{m>0} \zeta_{t,m} x^m, \qquad 0 \leq x \leq 1,\ 0 \leq t < +\infty.$$

Observe that due to the gelation effect at times $t>1$ the total mass $\zeta_{t,*} = u(1,t)$ is less than 1, hence $u$ is for $t>1$ no longer a probability generating function in the strict sense. The evolution of $u$ comprises all necessary information about the limit system including the mass of the gel-particle $\zeta_{t,\infty} = 1 - u(1,t)$ and the expected particle size per monomer $\sigma_t = \frac{\partial}{\partial x} u(1,t)$. From a well-known formula for the probability generating function of the total progeny of a Galton–Watson process (cf. [8]) we derive the following semi-explicit formula for $u$:

$$x = u e^{t(1-u)}. \tag{11}$$

Now observe some interesting consequences of this solution. First observe that for each $t$ (11) has the trivial solution $u=1$ for $x=1$, and it is not difficult to check that for $t \leq 1$ this is the only solution of $1 = u e^{t(1-u)}$ fulfilling $0 \leq u \leq 1$, while for $t>1$ there is a second solution with $0 < u < 1$ representing the true value of $u(1,t) = \zeta_{t,*} = 1 - \zeta_{t,\infty}$; see Figures 2 and 3.

Observe that the singularity at $t=1$ appears as a cusp catastrophe in the full solution of (11), and the existence of the gel-particle seems to find some reflection in the solution $u(x,t) \equiv 1$.

So for $t \geq 1$ the total mass is given by the inverse of the function

$$u \to t = \frac{-\log u}{1-u}.$$

Also, it is easy to derive from this that the initial slope of $\zeta_{t,*}$ at $t=1$ is $-2$, respectively the initial growth speed of the gel-particle is 2, a fact mentioned in [2], page 123 (where it has to be taken into account that the time scaling there differs from ours by a factor of 2, such that his value is 4). Then, in course of time $\zeta_{t,*}$ tends to zero exponentially fast.

Another implication of (11) is a relation for the expected particle size per monomer $\sigma_t$ which reads as $\sigma_t = \frac{\zeta_{t,*}}{1 - t \zeta_{t,*}}$. From this it is easy to derive the asymptotics of $\sigma$ in a right neighborhood of 1 as $\sigma_t \sim (t-1)^{-1}$ which shows that $|t-1|^{-1}$ is the asymptotics for $\sigma$ near $t=1$ (for $t<1$ this formula is well known to be the exact expression for $\sigma$).

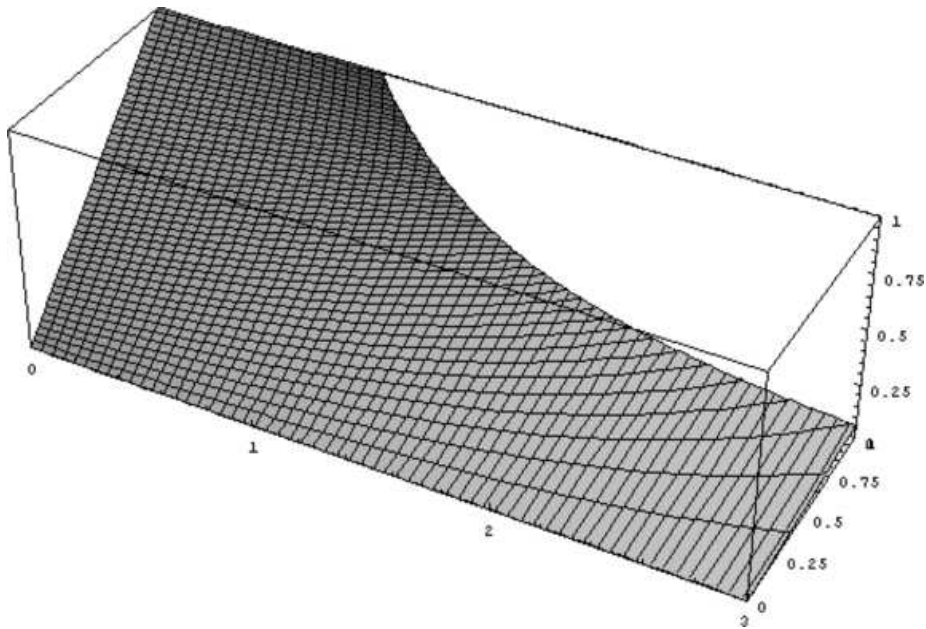

Fig. 2. *Plot of $u(x,t)$, $0 \leq t \leq 3, 0 \leq x \leq 1$.*

**4. Expected behavior of the spatial model after induced gelation.** The following remarks are only heuristical, but they are supported by the results of numerical simulations.

It is clear that *after* gelation the spatial model will behave in a more complex way, since the limiting model cannot be expected to be nonrandom. In fact, there will be two gel-particles at the two sites immediately after gelation, and these are not subject to a law of large numbers in their migration behavior. Consequently, the limiting behavior after $t_{\text{gel}}$ must be expected to be characterized by jumps of the giant particle(s) from one site to the other. They immediately coalesce if they happen to be at the same site. The further time evolution of the quantities $c_{t,m,j}^{\lambda}$ strongly depends on the presence or absence of a giant particle at $j$.

Hence, in contrast to the (mean-field) stochastic coalescent the (asymptotic) gelation effect does not only imply a *loss of total mass* (caught by the giant particle) but additionally the spatial model shows up a *stochastic limiting behavior*. So the limit model cannot be described by a (deterministic) kinetic equation similar to (8).

In order to describe the limiting object as a Markov process it is convenient to include the two gel-particles with their masses. This is not necessary in the mean-field model since the gel-particle simply represents the lost mass, whereas in our situation the distribution of the lost mass between the giant particles and their random locations are important.

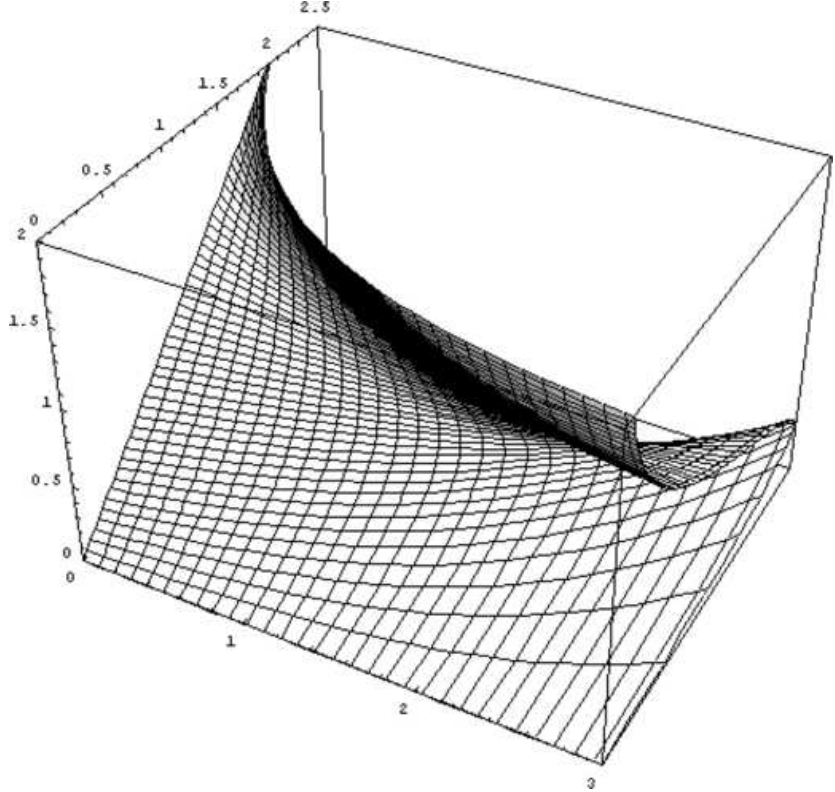

Fig. 3. *Extended plot of $u(x,t)$, $0 \leq t \leq 3, 0 \leq x \leq 2.5$.*

In order to unify the setup, we consider from this moment rescaled masses $N^{-1} m \xi_{t,m,j}$ instead of rescaled particle numbers $N^{-1} \xi_{t,m,j}$. With this change of the point of view the Smoluchowski system valid for $t \leq t_{\text{gel}}$ gets modified to a mass flow equation (writing $\zeta_{t,m,j}^\lambda$ for $m c_{t,m,j}^\lambda$)

$$\dot{\zeta}_{t,m,j}^\lambda = \tfrac{1}{2} m \sum_{0<n<m} \zeta_{t,n,j}^\lambda \zeta_{t,m-n,j}^\lambda - m \zeta_{t,m,j}^\lambda \sum_{n>0} \zeta_{t,n,j}^\lambda - \kappa \zeta_{t,m,j}^\lambda + \kappa \zeta_{t,m,1-j}^\lambda,$$

for $0 \leq t \leq t_{\text{gel}}$. Let us denote the masses of the gel-particles by $\zeta_{t,\infty,j}^\lambda, i \in \{0,1\}$. They are zero for $t \leq t_{\text{gel}}$. Finally, let $\pi_j(t)$, $j \in \{0,1\}$, $t \geq 0$ be two completely independent renewal processes with exponential waiting time (parameter $\kappa^{-1}$), with $\pi_j(t)$ denoting the time since the last renewal event occurred. These two processes are the only random source of the limiting model. We get the stochastic evolution equations

$$\dot{\zeta}_{t,m,j}^\lambda = \tfrac{1}{2} m \sum_{0<n<m} \zeta_{t,n,j}^\lambda \zeta_{t,m-n,j}^\lambda \tag{12}$$

$$- m\zeta_{t,m,j}^\lambda \left( \zeta_{t,\infty,j}^\lambda + \sum_{n>0} \zeta_{t,n,j}^\lambda \right) - \kappa \zeta_{t,m,j}^\lambda + \kappa \zeta_{t,m,1-j}^\lambda$$

and

$$\begin{aligned}
\zeta_{t,\infty,j}^\lambda = \mathbf{1}_{\{\pi_{1-j}(t) \geq \pi_j(t)\}} \int_{t-\pi_j(t)}^t \zeta_{s,\infty,j}^\lambda \sum_{n>0} n \zeta_{s,n,j}^\lambda \, ds \\
+ \mathbf{1}_{\{\pi_{1-j}(t) < \pi_j(t)\}} \Bigg[ \zeta_{t-\pi_{1-j}(t),\infty,j}^\lambda + \zeta_{t-\pi_{1-j}(t),\infty,1-j}^\lambda \\
+ \int_{t-\pi_{1-j}(t)}^t \zeta_{s,\infty,j}^\lambda \sum_{n>0} n \zeta_{s,n,j}^\lambda \, ds \Bigg],
\end{aligned}$$
(13)

for $t \geq 0$. This means in more rigorous terms that the limiting process is a Markov jump process with piecewise continuous trajectories, where the jumps are determined by the (external) renewal processes $\pi_j(t)$, and between jumps the trajectory is described by an ODE (i.e., deterministic), while at the discrete jump times $t_1, t_2, \ldots$ given by $\pi_j(t)$ a gel-particle of mass $\zeta_{t_i,\infty,j}^\lambda$ changes from $j$ to $1-j$ and forms at site $1-j$ together with the gel-particle there (if there is one) a larger gel-particle of mass $\zeta_{t_i,\infty,j}^\lambda + \zeta_{t_i,\infty,1-j}^\lambda$. Observe that we need to write the deterministic evolution equation between jumps as an integral equation since another gelation may have time to occur at the given site as long as $\zeta_{t,\infty,j}^\lambda$ is zero with the consequence that discontinuities in the derivative occur. The trajectories of $\zeta_{t,m,j}^\lambda$ are piecewise differentiable between the jumps.

So, from gelation time up to the first jump time happening after gelation [the pre-gelation jumps given by $\pi_j(t)$ are irrelevant] the whole evolution is completely deterministic and may be described by (12) plus the equation for the gel-particle's masses (13) which can be rewritten as $\dot{\zeta}_{t,\infty,j}^\lambda = \zeta_{t,\infty,j}^\lambda \sigma_{t,j}^\lambda$.

Just as in the mean-field model *after* gelation equation (4) has to be modified, taking into account the strong effect of the gel-particle, and ceases to be autonomous:

$$\dot{\sigma}_{t,j}^\lambda = (\sigma_{t,j}^\lambda)^2 + \kappa(\sigma_{t,1-j}^\lambda - \sigma_{t,j}^\lambda) - \zeta_{t,\infty,j}^\lambda \sum_{n>0} n^2 \zeta_{s,n,j}^\lambda.$$

It involves the second moment of the mass distribution, respectively the "third moment" in the traditional notation (referring to $c$ instead of $\zeta$). So, while $\sigma \zeta_\infty$ gives the speed of growth of $\zeta_\infty$ (of the gel-particle), the speed of descent of $\sigma$ now includes $\rho \zeta_\infty$, where $\rho$ is the second moment of the mass distribution, this being valid in the mean-field model as well as in our two-site model. The above differential equation for $\sigma$ may be expected to be valid except for the jump times (where right and left derivative are different) *and* except for the gelation time(s) (in our case there is more than one of them with positive probability as explained above). In the latter case it becomes meaningless since $\sigma$ exhibits a singularity (and the indefinite expression $0 \cdot \infty$ is involved).

In view of the fact that the consideration of the (probability) generating function of the mass distribution in the one-site (homogeneous) situation gives a detailed picture of the behavior of the limit model also after gelation, it is tempting to try this method also in our case to retrieve information about the deterministic behavior until the first jump happens. So one might hope for some more information about the asymptotics of $\zeta_\infty$ and $\sigma$ after $t_{\rm gel}$. Of course, it should clearly be expected that at the *inducing site* the asymptotics of $\sigma$ remains unchanged as in the one-site case $|t - t_{\rm gel}|^{-1}$ due to the fact that the influence of the *induced site* is negligible near $t_{\rm gel}$. It is less clear what will happen at the induced site in a *right* neighborhood of $t_{\rm gel}$, since here the effects of feeding $\sigma$ (at a logarithmic order, by the still large value of its counterpart at the inducing site) and of reducing $\sigma$ by the (comparatively small) gel-particle are competitive.

Let us introduce the generating functions

$$u(x,t,j) = \sum_m \zeta^\lambda_{t,m,j} x^m.$$

Then the evolution equations for $\zeta^\lambda_{t,m,j}$ imply the following equation for $u$, which should be valid until the first jump happens (after that $u$ itself will be random):

$$\frac{\partial}{\partial t} u(x,t,j) = x(u(x,t,j) - \zeta^\lambda_{t,\infty,j} - u(1,t,j)) \frac{\partial}{\partial x} u(x,t,j)$$
$$+ \kappa(u(x,t,1-j) - u(x,t,j)).$$

Observe that the first, nonlinear term differs from the one-site situation by the appearance of the total mass $\zeta^\lambda_{t,\infty,j} + u(1,t,j)$ at the given site instead of simply 1 in the one-site expression $x(u-1)\frac{\partial}{\partial x} u$. The new linear term is due to the motion of particles. Of course this equation is not autonomous and has to be supplemented by

$$\frac{d}{dt} \zeta^\lambda_{t,\infty,j} = \zeta^\lambda_{t,\infty,j} \cdot \frac{\partial}{\partial x} u(1,t,j),$$

which is just the rewritten evolution equation of $\zeta^\lambda_{t,\infty,j}$ from above. Obviously, this system of PDEs looks much less promising for having a (semi-)explicit solution, as was the case in the one-site model.

We emphasize again that all these considerations are heuristical though likely to be valid, until a rigorous proof of the convergence to the limiting process is available.

We defined gelation for our model in a weak sense meaning that $\sigma$ tends to infinity, which is strong evidence (but of course no rigorous proof) for the simultaneous appearance of gel-particles at the *first* gelation instance if $\kappa > 0$. In fact, the possible later gelation events cannot be induced gelations, since the presence of a gel-particle at the other site keeps the second moment there at finite values. So the first gelation is particularly interesting (and the proof of the rigorous limit theorem can be expected to be difficult mainly because of this special circumstance).

Bear in mind that, for $N$ being large but finite, only a tiny portion of logarithmic order of the quantity $\sigma$ feeds the induced gelation at the other site, or, to say it in different words, there is only a small time window around gelation time where many very large particles aggregate at the site where the "genuine" gelation just started, so only a small part of them can be expected to migrate to the other site to induce gelation there. In fact, it is quite likely that just as in the mean-field model the largest particles (except for the gel-particle) are of size $\leq N^{2/3}$ and these "mesoscopic" particles are present only for a short time before and after gelation. They carry the growth of the $\sigma$. The respective part of them which succeeds in jumping to the other site during their short time of existence is likely to be of logarithmic order in view of the calculations made in the proof of Theorem 4. This leads to the situation that, in contrast to the inducing site's growth of the expected particle mass per monomer to "infinity" (in fact, some power $N^\alpha$, $\alpha \leq \frac{2}{3}$), the growth of the expected particle mass at the induced site is carried only by the few very large particles which immigrated and will probably be of logarithmic order in $N$. Therefore three effects should be likely to be observed even for *very* large $N$ of the order of Loschmidt's constant or even of the order $10^{100}$:

(1) the induced gelation happens with a significant *delay* of the order $(\log N)^{-1}$;

(2) since it is carried by only few large "foreign" particles, in contrast to the inducing gelation the growth of $\sigma$ in the (delay) pre-gelation interval can be expected to show visible erratic deviations (caused by jumps from/to the other site and coalescences between the large particles) from the law of large numbers; and

(3) since the initial growth speed of the gel-particle is governed by the value of $\sigma$ at gelation time, which is maximally of logarithmic order in $N$ at the induced site, the second gel-particle's growth will be slow [we referred in Section 3 to the fact that in the mean-field case and hence probably also at the inducing site the (unrescaled) primary slope of the gel-particle's mass is of order $2N$].

Though at the time being no *quantum computer* was available to us to study simulations with $10^{100}$ particles (or only $10^{23}$), the simulations with $N = 10^9$ rather clearly show the predicted peculiarities of the induced gelation; see Figure 4 below. In this simulation, the behavior at site 0 is clearly governed by the law of large numbers, except for the random jump of the giant particle from 1 to 0 at time 5.86. At site 1 the boost of the second moment induced by site 0 around the (inducing) gelation time is visible, but not too strong. The delay between inducing and induced gelation is about one time unit, and during this delay time the second moment at site 1 is mainly carried by just a few particles, which fact causes strong fluctuations visible as some large jumps, respectively a "Brownian-like" behavior around the (induced) gelation. As criterion for gelation a threshold of size $2 \cdot 10^6$ for the largest particle was chosen, that is, $2N^{2/3}$.

**5. Remarks concerning other gelling kernels.** The situation for general gelling kernels (instead of $mn$) is more complicated, even in the mean-field case. Existence of gelling solutions has been conjectured for a long time (see the discussion in [1]). Rigorous results were given in [6, 10]. The asymptotic behavior of the stochastic model was studied, for example, in [3, 7, 10, 14, 15]. Convergence has been proved in the pre-gelation interval. After gelation, the Smoluchowski equation might have to be modified (dependent on the kernel). Convergence is established under the assumption that uniqueness holds. The asymptotic size of the largest particle is still an open problem. It has been studied in [1, 16] (see also the corresponding sections in [3, 19]).

Here we add some heuristic comments concerning the particle behavior around gelation, explaining the peculiarity of the $mn$-kernel in comparison to several other gelling kernels.

*Kernel $(mn)^\alpha$ for $1/2 < \alpha < 1$.* If there was a giant particle of size $\gamma N$, it would absorb the particles of fixed size $m$ at a rate

$$N^{-1}(\gamma m N)^\alpha = \frac{\gamma^\alpha m^\alpha}{N^{1-\alpha}},$$

that is, at an asymptotically vanishing rate. So it looks a bit surprising that it occurs at all. The explanation should be some kind of *food chain*. Just as in the $\alpha = 1$ case, in the pre-gelation phase the coagulation process shifts the total mass progressively to particles of larger size (individual mass), letting at the gelation point the mass flow "into infinity" reach an asymptotically positive speed. That means that the speed at which mass is lost from the set of all particles of size $m \leq B$ no longer goes to zero as $B$ tends to infinity (at a suitable slower speed than $N$). Hence a positive part of the mass accumulates at particles of a size which grows to infinity as $N$ goes to infinity. An indicator of that situation is that the mean particle size per monomer ("second moment" in the usual terminology) $\sigma$ grows to infinity. It is not a priori obvious what should be the typical size of the particles carrying that lost mass.

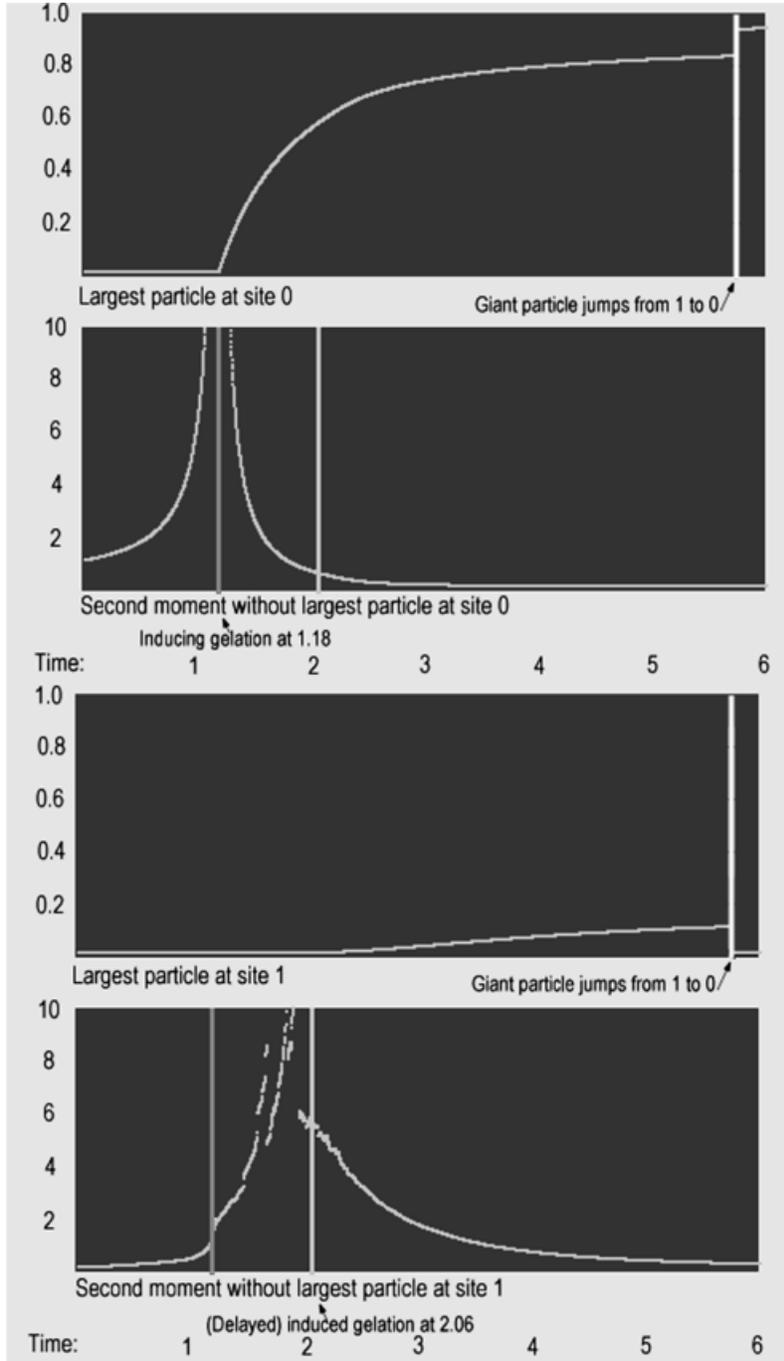

Fig. 4. *Simulation with one billion particles, $\lambda = 0, \kappa = 0.29$.*

Anyway, it might be conjectured that just as in the $\alpha = 1$ case the "lost mass" changes the regime abruptly by bringing the second moment back to finite values and stopping the "active" transfer (from the point of view of bounded-size particles) of mass to infinity. If this were true, the lost mass would actively absorb mass from the finite-size particles at positive speed, preferring large particles (that is the reason why the second moment drops to finite values for $\alpha = 1$; it should be noted that here "finite" always means: remaining bounded for $N$ tending to infinity, and does not refer to a limit model). As was mentioned above, the absorbing power of the giant particle is not enough to do this. For $\alpha < 1$ it must be taken into account that, while the rate at which a single particle of size $\gamma N$ absorbs $m$-particles is given by $\frac{\gamma^\alpha m^\alpha}{N^{1-\alpha}}$, for the same mass $\gamma N$ being distributed on $\gamma N^{1-\beta}$ particles of size $N^\beta$ the absorbing force is given by

$$N^{-1}\gamma N^{1-\beta}(mN^\beta)^\alpha = \gamma m^\alpha N^{\beta(\alpha-1)}.$$

The exponent is still negative, indicating that, though splitting the mass into many smaller particles increases the absorbing force, the totality of "infinite-size" particles is unable to actively withdraw mass from the bounded region.

Consequently, gelation will go on in a way very different from the critical case $\alpha = 1$: The bounded-size particles do actively transfer mass to infinity even after the gelation point, and this mass is transported through a kind of *food chain* at positive speed to infinity, reaching the giant particle only at the end of this chain. Observe that the mass transfer from particles of size $N^\beta$ directly to the giant size $N$ particle is given by a rate

$$N^{-1}(N^\beta N)^\alpha = N^{\alpha(1+\beta)-1}$$

which tends to zero unless $\beta \geq \frac{1}{\alpha} - 1$. So the food chain should be located between exponents 0 and $\frac{1}{\alpha} - 1$.

It should be expected that, in contrast to the critical case $\alpha = 1$, the second moment $\sigma$ keeps being infinite after the gelation point.

As for a corresponding spatial model (with mass-independent jump rate), the behavior should be significantly different compared to the critical $mn$ case: The portion of mass concentrated on finite-size particles might behave deterministically after gelation, too, because of a law of large numbers. This should be a kind of deterministic spatial Smoluchowski equation. Since the giant particle(s) at the spatial sites is/are unable to directly consume mass from below the threshold $\frac{1}{\alpha} - 1$, the random jumps performed by them (or it) should be without effect on the rest of the mass.

*Kernel $m^\alpha n + mn^\alpha$ for $\alpha < 1$.* The behavior will be similar to the $mn$-case, since the giant particle is able to withdraw mass from the bounded region at positive speed in this case.

*Kernel $m^\alpha + n^\alpha$ for $\alpha > 1$.* This kernel can be expected to be instantaneously and completely gelling [1, 11]. In fact: For large $N$ after an arbitrary short time an arbitrary large particle of size $N'$ is present with arbitrary large probability (of course at the beginning $1 \ll N' \ll N$). This particle by the law of large numbers grows like the solution of the ODE $\dot{x} = x^\alpha$, starting at initial value $N'$. So it takes an arbitrary short time to consume the rest of the population. Hence no interesting spatial model should be expected.

*Kernel $m^\alpha n^\beta + m^\beta n^\alpha$ for $\alpha + \beta > 1$.* If one of the parameters is greater than 1, the argument of the case $m^\alpha + n^\alpha, \alpha > 1$ applies. Otherwise, if $\alpha, \beta < 1$, it seems that anything said about the $(mn)^\alpha$ case could be valid here, too.

**6. Proofs.** We start by quoting a result on differential inequalities in $\mathbb{R}^d$ which will be used in the following considerations at several places. Let for $x, y \in \mathbb{R}^d$ the notation $x \leq y$ mean that $x_i \leq y_i$ for $i = 1, 2, \ldots, d$. A function $g : \mathbb{R}^d \to \mathbb{R}^d$ is called *quasi-monotone increasing* (with respect to this natural semi-ordering) if for each $z \in \mathbb{R}^d$, $0 \leq z$ from $x \leq y$ and $\langle z, x \rangle = \langle z, y \rangle$ it follows that $\langle z, g(x) \rangle \leq \langle z, g(y) \rangle$. In fact, we shall employ this definition only for the cases $d = 1$ (where any function is quasi-monotone increasing) and $d = 2$, where it simply means that $g_i(x)$ is a monotone increasing function of $x_{3-i}$ when $x_i$ is fixed, for $i = 1, 2$. Then we have (cf. [9] or [17]):

THEOREM 5. *Let $f : [a, b] \times \mathbb{R}^d \to \mathbb{R}^d$ be a function such that $f(t, \cdot)$ is quasi-monotone increasing and locally Lipschitz continuous for each $t \in [a, b]$, and let $u, v : [a, b] \to \mathbb{R}^d$ be differentiable. Then*

$$\dot{u}(t) - f(t, u(t)) \leq \dot{v}(t) - f(t, v(t)), \qquad t \in [a, b]$$

*and*

$$u(a) \leq v(a)$$

*imply*

$$u(t) \leq v(t), \qquad t \in [a, b].$$

6.1. *Proof of Theorem 3.* First observe that the theorem about existence and uniqueness of ODE ensures that there is a unique maximal local smooth solution to this evolution equation up to some $0 < T_{\max} \leq +\infty$, and if $T_{\max} < +\infty$, then necessarily the vector $(x, y)$ becomes unbounded near $T_{\max}$. Consider $S := x + y$. It fulfills the equation $\dot{S} = x^2 + y^2 \geq \frac{1}{2}(x+y)^2 = \frac{1}{2}S^2$. The maximal solution of the corresponding ODE is $(1 - \frac{1}{2}t)^{-1}$, $0 \leq t < 2$. Hence we see (using Theorem 5 in the case $d = 1$) that $T_{\max} \leq 2$ is finite. Let $t_{\text{gel}} := T_{\max}$. We know that at least one of the $x, y$ will be unbounded near $t_{\text{gel}}$. Assume it is the function $x$. Since $\dot{x} \geq x^2 - \kappa x \geq -\kappa x$ and $\dot{y} \geq y^2 - \kappa y \geq -\kappa y$, the functions are nonnegative for $t > 0$ (again by Theorem 5). Furthermore, since $x$ is unbounded near $t_{\text{gel}}$, this estimate combined with Theorem 5 shows that it is in fact monotonically increasing near $t_{\text{gel}}$ and hence tends to $+\infty$. Assume $y$ is bounded near $t_{\text{gel}}$. Then there is some $\delta$ such that $y^2 + \kappa(x-y) > \frac{1}{2}\kappa x$ for $t > t_{\text{gel}} - \delta$ implying $\dot{y} > \frac{1}{2}\kappa x$ for $t > t_{\text{gel}} - \delta$.

Now represent $x$ as $g(t_{\text{gel}} - t)^{-1}$ in a left neighborhood of $t_{\text{gel}}$ with $g$ being smooth in a right neighborhood $(0, \delta)$ of $0$ where we may use the same $\delta$ as above. We get

$$\dot{x} = x^2 + \kappa(y - x) = g(t_{\text{gel}} - t)^{-2} \dot{g}(t_{\text{gel}} - t) = x^2 \dot{g}(t_{\text{gel}} - t),$$

so that

$$\dot{g}(t) = [(x(t_{\text{gel}} - t))^2 + \kappa(y(t_{\text{gel}} - t) - x(t_{\text{gel}} - t))]/(x(t_{\text{gel}} - t))^2 \to 1$$

as $t \to 0$ by assumption. Hence $g(t) = c(t)t$ with $c$ being bounded from above and below by a positive constant $C$ (resp. its inverse $C^{-1}$) in a right neighborhood of $0$. We derive

$$y(t_{\text{gel}} - t) = y(t_{\text{gel}} - \delta) + \int_{t_{\text{gel}} - \delta}^{t_{\text{gel}} - t} \dot{y}(u)\, du > y(t_{\text{gel}} - \delta) + \frac{\kappa}{2} \int_{t_{\text{gel}} - \delta}^{t_{\text{gel}} - t} x(u)\, du$$

$$= y(t_{\text{gel}} - \delta) + \frac{1}{2}\kappa \int_{t_{\text{gel}} - \delta}^{t_{\text{gel}} - t} c(t_{\text{gel}} - u)^{-1}(t_{\text{gel}} - u)^{-1}\, du$$

$$\geq \frac{1}{2}\kappa C \int_{t_{\text{gel}} - \delta}^{t_{\text{gel}} - t} (t_{\text{gel}} - u)^{-1}\, du = \frac{1}{2}\kappa C(\ln(\delta) - \ln t).$$

The last expression is unbounded as $t \to 0$ in contradiction to the assumption. Hence we know that $y$ is unbounded near $t_{\text{gel}}$, too. As above we may conclude that $y \to +\infty$ as $t \to t_{\text{gel}}$ monotonically. We have shown the assertion.

6.2. *Proof of Theorem* 4. Consider $D = x - y$, $S = x + y$ and assume $x_0 > y_0$. These functions fulfill the quadratic ODE

$$\dot{S} = \tfrac{1}{2}S^2 + \tfrac{1}{2}D^2, \qquad \dot{D} = SD - 2\kappa D$$

with initial condition $S_0 = 1$, $0 < D_0 \leq S_0$. Since $S$ is obviously monotonically increasing we get $\dot{D} = SD - 2\kappa D \geq (1 - 2\kappa)D$ from which we infer by Theorem 5 that $D > 0$ is always fulfilled. Hence the site with originally more monomers has always a larger expected particle size. And, since $S$ goes to infinity monotonically near $t_{\text{gel}}$, we conclude that $\dot{D}$ is positive near $t_{\text{gel}}$, so $D$ increases monotonically near $t_{\text{gel}}$, too. Assume $D$ is bounded near $t_{\text{gel}}$. As above, represent $S$ as $f(t_{\text{gel}} - t)^{-1}$ in a left neighborhood of $t_{\text{gel}}$ which yields

$$\dot{S} = \tfrac{1}{2}S^2 + \tfrac{1}{2}D^2 = f(t_{\text{gel}} - t)^{-2}\dot{f}(t_{\text{gel}} - t) = S^2\dot{f}(t_{\text{gel}} - t)$$

and we may infer again that $\dot f$ would tend to $\frac{1}{2}$ near 0 and as above we would obtain a contradiction, since in this case the integral

$$\int_{t_{\text{gel}}-\delta}^{t_{\text{gel}}-t} S(u)D(u)\,du \geq D(t_{\text{gel}}-\delta) \int_{t_{\text{gel}}-\delta}^{t_{\text{gel}}-t} S(u)\,du$$

diverges near $t_{\text{gel}}$. So in fact $D$ increases monotonically to infinity near $t_{\text{gel}}$.

Let $Q := \frac{D}{S} \leq 1$. Then it is easy to see that $Q$ fulfills

(14) $$\dot Q = \tfrac{1}{2}D(1-Q^2) - 2\kappa Q.$$

If $Q$ was bounded away from 1 near $t_{\text{gel}}$ we might conclude from this differential equation that it is strictly increasing near $t_{\text{gel}}$, hence having a limit $Q_0 < 1$. With

$$\dot S = \tfrac{1}{2}S^2 + \tfrac{1}{2}D^2 = f(t_{\text{gel}}-t)^{-2}\dot f(t_{\text{gel}}-t) = S^2\dot f(t_{\text{gel}}-t)$$

we would conclude that $\dot f$ would tend to $\frac{1}{2}(1+Q_0^2)$ near 0 and hence $\int_{t_{\text{gel}}-\delta}^{t_{\text{gel}}-t} S(u)\,du$ as well as $\int_{t_{\text{gel}}-\delta}^{t_{\text{gel}}-t} D(u)\,du$ would diverge near $t_{\text{gel}}$. This would imply the divergence of $Q$ near $t_{\text{gel}}$ in view of (14), a contradiction. We have shown that 1 is an accumulation point of $Q$ near $t_{\text{gel}}$. If it was not the limit, $\dot Q$ would necessarily assume arbitrary large negative values in contradiction to the boundedness from below by $-2\kappa Q$. We have shown that $Q \to 1$ as $t \to t_{\text{gel}}$. This means that $y/x \to 0$ as $t \to t_{\text{gel}}$, provided $x_0 > y_0$. Repeating the calculations some lines above, this time with the confirmed value $Q_0 = 1$, we get that $S = f(t_{\text{gel}}-t)^{-1}$, where $\dot f \to 1$ as $t \to 0$. So there are smooth (in a right neighborhood of—and excluding—0) functions $c, c_1, c_2$ all tending to 1 as $t \to 0$ such that

$$S = c_1(t_{\text{gel}}-t)(t_{\text{gel}}-t)^{-1}, \qquad D = c_2(t_{\text{gel}}-t)(t_{\text{gel}}-t)^{-1},$$
$$x = c(t_{\text{gel}}-t)(t_{\text{gel}}-t)^{-1},$$

if $x_0 > y_0$. It remains to calculate the asymptotics of $y$.

We consider the quantity $q := y/x < 1$ which fulfills $\dot q = q^2(x-\kappa) - qx + \kappa$. We know that $q$ is positive but tends to zero near $t_{\text{gel}}$, and $qx = y$ tends to infinity. Hence for $t$ sufficiently close to $t_{\text{gel}}$ we have $\dot q < -(1-\delta)qx$ with arbitrary small $\delta$. So

$$[\log q]^\cdot < -(1-\delta)x < -(1-\delta)c(t_{\text{gel}}-t)(t_{\text{gel}}-t)^{-1} < -(1-2\delta)(t_{\text{gel}}-t)^{-1}$$

for $t$ sufficiently close to $t_{\text{gel}}$. Consequently we get $\log q < c_0 + (1-2\delta) \times \log(t_{\text{gel}}-t)$ (for some $c_0$) in a neighborhood of $t_{\text{gel}}$. So $q < e^{c_0}(t_{\text{gel}}-t)^{1-2\delta}$ and hence $y < e^{c_0}(t_{\text{gel}}-t)^{-2\delta}$. This proves that $y$ is not only $o(x)$ but even $o(x^\varepsilon)$ for each $\varepsilon > 0$.

We return to the original equation $\dot y = y^2 + \kappa(x-y)$. Now it is obvious that for $y_0 < x_0$ the increase of $y$ near $t_{\text{gel}}$ is (asymptotically) completely caused by the migration term $\kappa x$, while the coagulation term $y^2$ can be neglected. For any positive $\delta$ we find for $t$ sufficiently close to $t_{\text{gel}}$: $(1-\delta)\kappa x < \dot{y} < (1+\delta)\kappa x$ and the $x$-asymptotics yields

$$-(1-2\delta)\kappa \log(t_{\text{gel}} - t) < k_0 + y < -(1+2\delta)\kappa \log(t_{\text{gel}} - t)$$

in a sufficiently small neighborhood of $t_{\text{gel}}$, where $k_0$ is a constant depending on $\delta$. Therefore we finally deduced the logarithmic asymptotics of $y = -c_3(t_{\text{gel}} - t)\kappa \log(t_{\text{gel}} - t)$ with $c_3$ tending to 1 as $t \to 0$.

### 6.3. Proof of Theorem 2.

*Existence.* Fix some $T$, $0 < T < t_{\text{gel}}^\lambda$. Consider the following system of ODEs:

$$\dot{c}_{t,m,j}^\lambda = \tfrac{1}{2} \sum_{0<n<m} n(m-n) c_{t,n,j}^\lambda c_{t,m-n,j}^\lambda - m c_{t,m,j}^\lambda \zeta_{t,*,j}^\lambda$$

(15)
$$- \kappa c_{t,m,j}^\lambda + \kappa c_{t,m,1-j}^\lambda,$$

$$\dot{\zeta}_{t,*,j}^\lambda = -\kappa \zeta_{t,*,j}^\lambda + \kappa \zeta_{t,*,1-j}^\lambda,$$

with initial condition

$$c_{0,1,0}^\lambda = \zeta_{0,*,0}^\lambda = \frac{1}{1+\lambda}, \qquad c_{0,1,1}^\lambda = \zeta_{0,*,1}^\lambda = \frac{\lambda}{1+\lambda}, \qquad c_{0,m,j}^\lambda = 0, \qquad m > 1.$$

It has obviously a unique solution for all $t$ since, given $B \in \mathbb{N}$, the equations for $m < B$, together with the equations for $\zeta_{t,*,j}^\lambda$ form a finite autonomous system. First observe that all quantities $c_{t,m,j}^\lambda, \zeta_{t,*,j}^\lambda$ are nonnegative. In fact, this is easily seen for $\zeta_{t,*,j}^\lambda$ which can be given in explicit terms:

(16) $\quad \zeta_{t,*,0}^\lambda = \frac{1}{2}\left(1 + \frac{1-\lambda}{1+\lambda}e^{-2\kappa t}\right), \qquad \zeta_{t,*,1}^\lambda = \frac{1}{2}\left(1 - \frac{1-\lambda}{1+\lambda}e^{-2\kappa t}\right).$

To see it for the $c_{t,m,j}^\lambda$, let us first modify the initial condition to $c_{0,m,j}^\lambda = c, m > 1$, where $c > 0$ is an arbitrary fixed value. It is enough to prove the nonnegativity of the corresponding solutions for all times $t > 0$ in these cases in view of the theorem about the continuous dependence of solutions of ODEs on the initial data. Assume the opposite. Then there would be some $c$ such that at least one of the $c_{t,m,j}^\lambda$ reaches 0 for the first time at some time $0 < t_0$. We get the differential inequalities $\dot{c}_{t,m,j}^\lambda \geq (-m - \kappa) c_{t,m,j}^\lambda$ being valid for $0 \leq t \leq t_0$ taking into account (16). But now we immediately infer from Theorem 5 that $c_{t_0,m,j}^\lambda > 0$ in contradiction to our assumption. Hence indeed all $c_{t,m,j}^\lambda$ are nonnegative for all times.

Next, consider the quantities

$$\zeta_{t,*,j}^{B,\lambda} := \sum_{m<B} m c_{t,m,j}^\lambda \quad \text{and} \quad \delta_{t,j}^{B,\lambda} := \zeta_{t,*,j}^{B,\lambda} - \zeta_{t,*,j}^\lambda.$$

We have

$$\dot{\zeta}_{t,*,j}^{B,\lambda} = \tfrac{1}{2} \sum_{k+l<B} (k+l)kl c_{t,k,j}^{\lambda} c_{t,l,j}^{\lambda} - \zeta_{t,*,j}^{\lambda} \sum_{m<B} m^2 c_{t,m,j}^{\lambda}$$

(17)
$$+ \kappa(\zeta_{t,*,1-j}^{B,\lambda} - \zeta_{t,*,j}^{B,\lambda})$$

$$\leq \sum_{k,l<B} k^2 l c_{t,k,j}^{\lambda} c_{t,l,j}^{\lambda} - \zeta_{t,*,j}^{\lambda} \sum_{m<B} m^2 c_{t,m,j}^{\lambda} + \kappa(\zeta_{t,*,1-j}^{B,\lambda} - \zeta_{t,*,j}^{B,\lambda})$$

$$= (\zeta_{t,*,j}^{B,\lambda} - \zeta_{t,*,j}^{\lambda}) \sum_{m<B} m^2 c_{t,m,j}^{\lambda} + \kappa(\zeta_{t,*,1-j}^{B,\lambda} - \zeta_{t,*,j}^{B,\lambda})$$

which implies

$$\dot{\delta}_{t,j}^{B,\lambda} \leq \delta_{t,j}^{B,\lambda} \sum_{m<B} m^2 c_{t,m,j}^{\lambda} + \kappa(\delta_{t,1-j}^{B,\lambda} - \delta_{t,j}^{B,\lambda}).$$

The initial values $\delta_{0,j}^{B,\lambda}$ are zero. So we may conclude, again by Theorem 5, that $\delta_{t,j}^{B,\lambda}$ is nonpositive yielding

(18)
$$\zeta_{t,*,j}^{B,\lambda} \leq \zeta_{t,*,j}^{\lambda} \leq 1.$$

Next, consider the quantity

$$\sigma_{t,j}^{B,\lambda} := \sum_{m<B} m^2 c_{t,m,j}^{\lambda}.$$

We have

$$\dot{\sigma}_{t,j}^{B,\lambda} = \tfrac{1}{2} \sum_{k+l<B} (k+l)^2 kl c_{t,k,j}^{\lambda} c_{t,l,j}^{\lambda} - \zeta_{t,*,j}^{\lambda} \sum_{m<B} m^3 c_{t,m,j}^{\lambda}$$

$$+ \kappa(\sigma_{t,1-j}^{B,\lambda} - \sigma_{t,j}^{B,\lambda})$$

(19)
$$\leq \sum_{k,l<B} k^3 l c_{t,k,j}^{\lambda} c_{t,l,j}^{\lambda} + (\sigma_{t,j}^{B,\lambda})^2 - \zeta_{t,*,j}^{\lambda} \sum_{m<B} m^3 c_{t,m,j}^{\lambda}$$

$$+ \kappa(\sigma_{t,1-j}^{B,\lambda} - \sigma_{t,j}^{B,\lambda})$$

$$= (\sigma_{t,j}^{B,\lambda})^2 + (\zeta_{t,*,j}^{B,\lambda} - \zeta_{t,*,j}^{\lambda}) \sum_{m<B} m^3 c_{t,m,j}^{\lambda} + \kappa(\sigma_{t,1-j}^{B,\lambda} - \sigma_{t,j}^{B,\lambda}).$$

Hence, by the above estimate $\zeta_{t,*,j}^{B,\lambda} \leq \zeta_{t,*,j}^{\lambda}$ we arrive at

$$\dot{\sigma}_{t,j}^{B,\lambda} \leq (\sigma_{t,j}^{B,\lambda})^2 + \kappa(\sigma_{t,1-j}^{B,\lambda} - \sigma_{t,j}^{B,\lambda}).$$

The function defined by $f_j(x) := x_j^2 + \kappa(x_{3-j} - x_j)$ is quasi-monotone increasing. Consequently by Theorem 5, $(\sigma_{t,j}^{B,\lambda})_{j=0,1}$ is dominated in $[0, t_{\text{gel}}^{\lambda})$ by the solution of the ODE (5)–(6), for all $B$. So $(\sigma_{t,j}^{\lambda})_{j=0,1}$ is bounded from above by this solution, too. Hence $\sigma_{t,j}^\lambda$ is integrable in any compact subinterval of $[0, t_{\text{gel}}^\lambda)$.

Next, consider again the relation (17) which implies

$$\zeta_{t,*,j}^{B,\lambda} = \zeta_{0,*,j}^{B,\lambda} + \int_0^t \left\{ \frac{1}{2} \sum_{k+l<B} (k+l) kl c_{s,k,j}^\lambda c_{s,l,j}^\lambda \right.$$

$$\left. - \zeta_{s,*,j}^\lambda \sum_{m<B} m^2 c_{s,m,j}^\lambda + \kappa(\zeta_{s,*,1-j}^{B,\lambda} - \zeta_{s,*,j}^{B,\lambda}) \right\} ds.$$

The quantity $\frac{1}{2}\sum_{k+l<B}(k+l)kl c_{s,k,j}^\lambda c_{s,l,j}^\lambda$ converges monotonically (as $B \to \infty$) toward

$$\tfrac{1}{2}\sum_{k,l}(k+l)kl c_{s,k,j}^\lambda c_{s,l,j}^\lambda = \sigma_{t,j}^\lambda \sum_m mc_{s,m,j}^\lambda$$

for each $s \in [0, t_{\text{gel}}^\lambda)$, in the same way $\sum_{m<B} m^2 c_{s,m,j}^\lambda$ converges toward $\sigma_{t,j}^\lambda$, and $\zeta_{s,*,j}^{B,\lambda}$ toward

$$\zeta_{s,*,j}^{\infty,\lambda} := \sum_m m c_{s,m,j}^\lambda \leq \zeta_{t,*,j}^\lambda \leq 1$$

[by (18)]. So from Lebesgue's theorem about monotone convergence we infer, taking into account the uniform boundedness of all these quantities in $[0, T]$, $T < t_{\text{gel}}^\lambda$,

(20) $\quad \zeta_{t,*,j}^{\infty,\lambda} = \zeta_{0,*,j}^{\infty,\lambda} + \int_0^t \{\sigma_{s,j}^\lambda(\zeta_{s,*,j}^{\infty,\lambda} - \zeta_{s,*,j}^\lambda) + \kappa(\zeta_{s,*,1-j}^{\infty,\lambda} - \zeta_{s,*,j}^{\infty,\lambda})\} ds.$

Let $\delta_{s,j}^{\infty,\lambda} := \zeta_{s,*,j}^{\infty,\lambda} - \zeta_{s,*,j}^\lambda \leq 0$ [by (18)] and $d_s^\lambda = |\delta_{s,0}^{\infty,\lambda} + \delta_{s,1}^{\infty,\lambda}| = -\delta_{s,0}^{\infty,\lambda} - \delta_{s,1}^{\infty,\lambda}$ and we get [bearing in mind that for each $s$ we have $\zeta_{s,*,0}^\lambda + \zeta_{s,*,1}^\lambda = 1$ in view of (15)]

(21) $\qquad d_t^\lambda = \int_0^t \{-\sigma_{s,0}^\lambda \delta_{s,0}^{\infty,\lambda} - \sigma_{s,1}^\lambda \delta_{s,1}^{\infty,\lambda}\} ds \leq C \int_0^t d_s^\lambda \, ds$

for some constant $C$ depending on $T$, since $\sigma_{s,j}^\lambda$ is uniformly bounded for $0 \leq s \leq T < t_{\text{gel}}^\lambda$. Since $d_s^\lambda$ is bounded by 2, we get from (21) for $s \leq \frac{1}{2C}$ that $d_s^\lambda$ is bounded by 1 and by induction that it must be zero in $[0, \frac{1}{2C}]$, which, arguing by repetition, leads to the conclusion that $d_s^\lambda$ is zero everywhere in $[0, t_{\text{gel}}^\lambda)$. Hence we have $\zeta_{s,*,j}^{\infty,\lambda} \equiv \zeta_{s,*,j}^\lambda$ everywhere in $[0, t_{\text{gel}}^\lambda)$. Now, substituting $\zeta_{s,*,j}^{\infty,\lambda} = \sum_m mc_{s,m,j}^\lambda$ for $\zeta_{s,*,j}^\lambda$ in (15) we see that the unique solution of (15) also satisfies (2) in $[0, t_{\text{gel}}^\lambda)$. We proved that this solution is positive and strong in the sense of integrability of $\sigma_{t,j}^\lambda$ in any compact interval contained in $[0, t_{\text{gel}}^\lambda)$.

*Uniqueness.* Let a strong positive solution $(c^\lambda_{t,m,j})_{m\in\mathbb{N},j=0,1,t\in[0,\tilde{t})}, \tilde{t} < t_{\text{gel}}$, of (2) be given. Using again the notation $\zeta^{B,\lambda}_{t,*,j} := \sum_{m<B} mc^\lambda_{t,m,j}$ and $\zeta^{\infty,\lambda}_{t,*,j} := \sum_m mc^\lambda_{t,m,j}$, from (2) we infer in the same way as in (17)

$$\dot\zeta^{B,\lambda}_{t,*,j} \leq (\zeta^{B,\lambda}_{t,*,j} - \zeta^{\infty,\lambda}_{t,*,j}) \sum_{m<B} m^2 c^\lambda_{t,m,j} + \kappa(\zeta^{B,\lambda}_{t,*,1-j} - \zeta^{B,\lambda}_{t,*,j})$$

which in view of $\zeta^{B,\lambda}_{t,*,j} \leq \zeta^{\infty,\lambda}_{t,*,j}$ yields $\dot\zeta^{B,\lambda}_{t,*,j} \leq \kappa(\zeta^{B,\lambda}_{t,*,1-j} - \zeta^{B,\lambda}_{t,*,j})$ and consequently $\dot\zeta^{B,\lambda}_{t,*,0} + \dot\zeta^{B,\lambda}_{t,*,1} \leq 0$. Thus the sum of the two positive quantities $\zeta^{B,\lambda}_{t,*,j}, j = 0,1$ is always less than or equal to its initial value 1 (for $B>1$). Letting $B$ tend to infinity we arrive at the conclusion that both $\zeta^{\infty,\lambda}_{t,*,0}, \zeta^{\infty,\lambda}_{t,*,1}$ are not greater than 1 for $0 \leq t < \tilde{t}$. Now we have

$$\zeta^{B,\lambda}_{t,*,j} = \zeta^{B,\lambda}_{0,*,j} + \int_0^t \tfrac{1}{2} \sum_{k+l<B} (k+l)kl c^\lambda_{s,k,j} c^\lambda_{s,l,j}\, ds$$
$$- \int_0^t \zeta^{\infty,\lambda}_{s,*,j} \sum_{m<B} m^2 c^\lambda_{s,m,j}\, ds + \int_0^t \kappa(\zeta^{B,\lambda}_{s,*,1-j} - \zeta^{B,\lambda}_{s,*,j})\, ds.$$

With $B \to \infty$ the first integral on the right-hand side tends toward

$$\int_0^t \tfrac{1}{2} \sum_{k,l}(k+l)kl c^\lambda_{s,k,j} c^\lambda_{s,l,j}\, ds = \int_0^t \sigma^\lambda_{s,j} \zeta^{\infty,\lambda}_{s,*,j}\, ds \leq \int_0^t \sigma^\lambda_{s,j}\, ds < +\infty$$

for each $0 \leq t < \tilde{t}$ by the assumption that the solution is strong. The same argument yields the finiteness of the second integral, which tends to the same limit $\int_0^t \sigma^\lambda_{s,j} \zeta^{\infty,\lambda}_{s,*,j}\, ds$. So we arrive at the system of two equations

$$\zeta^{\infty,\lambda}_{t,*,j} = \zeta^{\infty,\lambda}_{0,*,j} + \int_0^t \kappa(\zeta^{\infty,\lambda}_{s,*,1-j} - \zeta^{\infty,\lambda}_{s,*,j})\, ds, \qquad j = 0,1,$$

from which we may easily conclude that the bounded quantities $\zeta^{\infty,\lambda}_{t,*,j}$ are continuous and even continuously differentiable and fulfill the ODE

$$\dot\zeta^{\infty,\lambda}_{t,*,j} = -\kappa \zeta^{\infty,\lambda}_{t,*,j} + \kappa \zeta^{\infty,\lambda}_{t,*,1-j}$$

with initial condition

$$\zeta^{\infty,\lambda}_{0,*,0} = \frac{1}{1+\lambda}, \qquad \zeta^{\infty,\lambda}_{0,*,1} = \frac{\lambda}{1+\lambda}.$$

Therefore our solution also fulfills (15) which has a unique solution for all $t$. So the maximal strong positive solution must be unique.

*Maximality.* We still have to show that the maximality implies $\tilde{t} = t_{\text{gel}}^\lambda$. Knowing already that the solution to (15) in $[0, t_{\text{gel}}^\lambda)$ is positive and strong [in view of the fact derived above that $\sigma_{t,j}^\lambda$ is bounded from above by the continuous solution of (4)], we just need to show that $\tilde{t}$ cannot be larger than $t_{\text{gel}}^\lambda$. For this sake we first prove that the second moments $(\sigma_{t,j}^\lambda)_{j=0,1}$, as announced above, fulfill (4) for $0 \leq t < t_{\text{gel}}^\lambda$. Denote by $\rho_{t,j}^{B,\lambda}, \rho_{t,j}^\lambda$ the quantities $\sum_{m<B} m^3 c_{t,m,j}^\lambda, \sum_m m^3 c_{t,m,j}^\lambda$, respectively, given the solution to (15). We get, in an analogous way as in (19) for $0 \leq t \leq T < t_{\text{gel}}^\lambda$,

$$\dot{\rho}_{t,j}^{B,\lambda} = \tfrac{1}{2} \sum_{k+l<B} (k+l)^3 k l c_{t,k,j}^\lambda c_{t,l,j}^\lambda - \zeta_{t,*,j}^\lambda \sum_{m<B} m^4 c_{t,m,j}^\lambda$$
$$+ \kappa(\rho_{t,1-j}^{B,\lambda} - \rho_{t,j}^{B,\lambda})$$
$$\leq \sum_{k,l<B} k^4 l c_{t,k,j}^\lambda c_{t,l,j}^\lambda + 3 \sum_{k,l<B} k^3 l^2 c_{t,k,j}^\lambda c_{t,l,j}^\lambda - \zeta_{t,*,j}^\lambda \sum_{m<B} m^4 c_{t,m,j}^\lambda$$
$$+ \kappa(\rho_{t,1-j}^{B,\lambda} - \rho_{t,j}^{B,\lambda})$$
$$= 3\rho_{t,j}^{B,\lambda} \sigma_{t,j}^{B,\lambda} + (\zeta_{t,*,j}^{B,\lambda} - \zeta_{t,*,j}^\lambda) \sum_{m<B} m^4 c_{t,m,j}^\lambda + \kappa(\rho_{t,1-j}^{B,\lambda} - \rho_{t,j}^{B,\lambda})$$
$$\leq 3\rho_{t,j}^{B,\lambda} C + \kappa(\rho_{t,1-j}^{B,\lambda} - \rho_{t,j}^{B,\lambda}),$$

where $C$ is a constant depending on $T$, where we made use of the boundedness of $\sigma_{t,j}^\lambda$ in compact subintervals of $[0, t_{\text{gel}}^\lambda)$. From this, in the same way as for the second moment, we derive that all $\rho_{t,j}^{B,\lambda}$ and hence also $\rho_{t,j}^\lambda$ are bounded by the solution of the coupled system of ODEs

$$\dot{x} = 3Cx + \kappa(y-x), \qquad \dot{y} = 3Cy + \kappa(x-y)$$

with $x_0 = \frac{1}{1+\lambda}, y_0 = \frac{\lambda}{1+\lambda}$ in $[0, T]$. Now we may rewrite the first equation in (19):

$$\sigma_{t,j}^{B,\lambda} = \sigma_{0,j}^{B,\lambda} + \int_0^t \tfrac{1}{2} \sum_{k+l<B} (k+l)^2 k l c_{s,k,j}^\lambda c_{s,l,j}^\lambda \, ds$$
$$- \int_0^t \zeta_{s,*,j}^\lambda \rho_{s,j}^{B,\lambda} \, ds + \kappa \int_0^t (\sigma_{s,1-j}^{B,\lambda} - \sigma_{s,j}^{B,\lambda}) \, ds.$$

With $B \to \infty$ the first integral on the right-hand side tends monotonically toward

$$\int_0^t \tfrac{1}{2} \sum_{k,l} (k+l)^2 k l c_{s,k,j}^\lambda c_{s,l,j}^\lambda \, ds = \int_0^t (\zeta_{s,*,j}^\lambda \rho_{s,j}^\lambda + (\sigma_{s,j}^{B,\lambda})^2) \, ds$$

being finite for $t < t_{\text{gel}}^\lambda$, and the second integral converges monotonically toward

$$\int_0^t \zeta_{s,*,j}^\lambda \rho_{s,j}^\lambda \, ds$$

which is finite for these $t$, too. So we may pass to the limit on both sides and get

$$\sigma_{t,j}^\lambda = \sigma_{0,j}^\lambda + \int_0^t (\sigma_{s,j}^\lambda)^2 \, ds + \kappa \int_0^t (\sigma_{s,1-j}^\lambda - \sigma_{s,j}^\lambda) \, ds.$$

Thus the $\sigma_{t,j}^\lambda$ are continuous and even continuously differentiable in $[0, t_{\text{gel}}^\lambda)$ and fulfill (4). Now from Theorem 4 we easily infer that the condition that $\sigma_{t,j}^\lambda$ be integrable for each compact subinterval of $[0, \tilde{t})$ (strong solution) implies $\tilde{t} \leq t_{\text{gel}}^\lambda$.

6.4. *Proof of Theorem* 1. We are now prepared to prove Theorem 1. At this point we take the point of view that the configuration at time $t$ of our particle system is given by a finite counting measure on $\mathbb{N} \times \{0, 1\}$,

(22) $$X_t = \sum_{m,j} \xi_{t,m,j} \delta_{(m,j)}.$$

Let us, just for the moment, renounce the assumption that the initial configuration be monomeric; that is, $\xi_{0,m,j}$ may be positive for $m > 1$.

In order to prove our limit theorem we consider a coupled system of truncated versions of the process $\{X_t\}$ indexed by $B \in \overline{\mathbb{N}} = \mathbb{N} \cup \{\infty\}$ where the dynamics of $\{X_t^B\}$ are given as follows:

1. First define

$$X_0^B := \int_{\mathbb{N} \times \{0,1\}} \mathbf{1}_{m<B} \delta_{(m,i)} X_0(d(m,i)) + \int_{\mathbb{N} \times \{0,1\}} \mathbf{1}_{m \geq B} m \delta_{(B,0)} X_0(d(m,i)),$$

that is, the particles of sizes less than $B$ remain unchanged, while at each site the masses of all particles larger than $B$ are combined and this mass is assigned to $(B, 0)$. $X_0$ can be represented as a finite sum $\sum_{l \in J} \delta_{(m_l, j_l)}$, where $J$ is a finite index set, and each $X_0^B$ can be written as

(23) $$\sum_{l \in J(B)} \delta_{(m_l, j_l)} + X_0^B(\{(B, 0)\}) \delta_{(B, 0)},$$

where $J(B)$ is the subset of $J$ corresponding to particles of size smaller than $B$. Obviously $X_0^\infty = X_0$.

2. It is easily seen that, by simply cutting each $X_t$ in an analogous way, the corresponding quantities will *not* form a Markov process except for $\kappa = 0$ (which corresponds to the one-site case). We intend to let $X_0^B$ evolve Markovian, thereby keeping it in close distance to $X_t$ for $B$ large, that is, to introduce a coupled evolution of $(X_t, X_t^B)$ where both marginal processes are Markovian and their distance can be controlled to be small for large $B$. To perform this, first introduce a new discrete parameter assigned to the particles taking two possible values v and h ("visible," "hidden"). That means, instead of $\mathbb{N} \times \{0,1\}$ we choose $\mathbb{N} \times \{0,1\} \times \{v, h\}$ as phase space. All particles present in $X_0$ with masses below $B$ are marked visible; the rest (with masses $\geq B$) is marked hidden. This gives $\widehat{X}_0^B$. Given a path $(X_t)_{t \geq 0}$ we construct $(\widehat{X}_t^B)_{t \geq 0}$ as follows, not changing the projection to $\mathbb{N} \times \{0,1\}$ which is $(X_t)_{t \geq 0}$:

(a) the mark "hidden" is dominant; that is, the combination of two particles occurring in $(X_t)_{t \geq 0}$ produces a hidden particle iff at least one of them was hidden, and jumping to the other site does not change the status hidden or visible,

(b) if a particle of mass at least $B$ is generated by the combination of two particles, it is marked hidden, and

(c) additionally, particles become "hidden" by a newly introduced interaction between particles of *different* sites: given $\widehat{X}_t^B$, a "visible" particle of mass $m < B$ at site $j$ becomes "hidden" at rate

$$m \int_{\mathbb{N} \times \{0,1\} \times \{v,h\}} k \mathbf{1}_{i=1-j, q=h} \widehat{X}_t^B(d(k,i,q)).$$

So the interaction concerns only the newly introduced mark, and has the effect that any pair of particles with masses $m, k$ *at different sites* assigns the mark "hidden" to the first particle at rate $mk$, assuming the second particle is already hidden. This represents a Markovian refinement of the original process.

Now we define $X_t^B$ by

$$X_t^B = \int_{\mathbb{N} \times \{0,1\} \times \{v,h\}} (\delta_{(k,i)} \mathbf{1}_{k<B, q=v} + \delta_{(B,0)} k \mathbf{1}_{q=h}) \widehat{X}_t^B(d(k,i,q)).$$

This means that only visible particles are present in $X_t^B$ at locations $< B$, whereas all hidden particles (in particular those with masses at least $B$) are collected to form the mass at $B$ at site 0. Observe that the total rate at which a visible particle of mass $m$ becomes "hidden" is given by the product of $m$ with the total mass of all other hidden particles.

3. So, if we do not intend to describe the joint dynamics of $(X_t, X_t^B)$ but only that of $X_t^B$ we may do it as follows, referring to the representation (23): Consider a family $(T_{l,l'}^B)_{l,l' \in J(B), l<l', j_l=j_{l'}}$ of independent exponentially distributed random variables with parameter $m_l m_{l'}$, a second (independent of the first) family $(S_l^B)_{l \in J(B)}$ of independent exponentially distributed random variables with parameter $\kappa N$, and a third independent family $(U_l^B)_{l \in J(B)}$ of independent exponentially distributed random variables with parameter $m_l X_0^B(\{(B,0)\})$. These families are finite in view of the finiteness of $J$. Let

$$T = \min_{l,l'}(T_{l,l'}^B, S_l^B, U_l^B).$$

Define $X_t^B = X_0^B$ for $t < T$ and

$$X_T^B = \begin{cases} X_0^B - \delta_{(m_l,j_l)} - \delta_{(m_{l'},j_{l'})} + \delta_{(m_l+m_{l'},j_l)}, & \text{if } T = T_{l,l'}^B, \\ & m_l + m_{l'} < B, \\ X_0^B - \delta_{(m_l,j_l)} - \delta_{(m_{l'},j_{l'})} + (m_l+m_{l'})\delta_{(B,0)}, & \text{if } T = T_{l,l'}^B, \\ & m_l, m_{l'} < B, \\ & m_l + m_{l'} \geq B, \\ X_0^B - \delta_{(m_l,j_l)} + \delta_{(m_l,1-j_l)}, & \text{if } T = S_l^B, m_l < B, \\ X_0^B - \delta_{(m_l,j_l)} + m_l\delta_{(B,0)}, & \text{if } T = U_l^B, \\ X_0^B, & \text{otherwise.} \end{cases}$$

Starting this construction anew with $X_T^B$ instead of $X_0^B$ we finally define $X_t^B$ for arbitrary $t \geq 0$, thus arriving at a Markov jump process $\{X_t^B\}$. Observe that $\{X_t^\infty\} = \{X_t^{\nu(X_0)+1}\} = \{X_t\}$, where $\nu(X_0) := \int_{\mathbb{N}\times\{0,1\}} m X_0(d(m,i))$ is the total mass.

The evolution of a single $X_t^B$ may be described as follows:

- Particles of sizes $m, n$ at the same site coalesce at rate $mn$, assuming $m, n < B$. If $m + n < B$, a new particle of size $m + n$ is generated, whereas for $m + n \geq B$ the mass $m + n$ is added to the mass already aggregated at $(B, 0)$.
- Particles of size less than $B$ change their location at rate $\kappa N$.
- Finally, particles of size $m < B$ are deleted at a rate given by the product of their size and the mass already aggregated at $(B, 0)$ and their mass is added to the aggregated mass.

The evolution conserves the total mass

$$\nu(X_0) \equiv \nu_B(X_0^B) := \int_{\mathbb{N}\times\{0,1\}, m<B} m X_0^B(d(m,i)) + X_0^B(\{(B,0)\}),$$

being independent of $B$.

Let, for $B \in \mathbb{N}$, $\mathcal{M}_B$ denote the space of finite signed measures on

$$(\mathbb{N} \cap [1, B]) \times \{0, 1\}$$

and define $L^B : \mathcal{M}_B \to \mathcal{M}_B$ by

$$\langle f, L^B(\mu) \rangle$$
$$= \tfrac{1}{2} \int_{(\mathbb{N}\times\{0,1\})^2} [f(x+y, j)\mathbf{1}_{x+y<B}$$
$$\qquad + f(B, 0)(x+y)\mathbf{1}_{x+y\geq B} - f(x, j) - f(y, j)]$$

$$\times \mathbf{1}_{x,y<B, j=k} \cdot xy(\mu \times \mu)(d[x,j],[y,k])$$

$$+ \int_{\mathbb{N}\times\{0,1\}} \mu(\{(B,0)\})[f(B,0)x - f(x,j)]\mathbf{1}_{x<B}x\mu(d[x,j])$$

$$+ \kappa \int_{\mathbb{N}\times\{0,1\}} [f(x,1-j) - f(x,j)]\mathbf{1}_{x<B}\mu(d[x,j]),$$

for any bounded function $f$ on $\mathbb{N} \times \{0,1\}$. Fix a finite measure $\mu_0$ on $\mathbb{N} \times \{0,1\}$ with

$$\nu(\mu_0) = \int_{\mathbb{N}\times\{0,1\}} x\mu_0(d[x,j]) < +\infty$$

and define

$$\mu_0^B := \int_{\mathbb{N}\times\{0,1\}} \mathbf{1}_{m<B}\delta_{(m,i)}\mu_0(d(m,i))$$

$$+ \int_{\mathbb{N}\times\{0,1\}} \mathbf{1}_{m\geq B}m\delta_{(B,0)}\mu_0(d(m,i)).$$

It is easy to see that for each $B \in \mathbb{N}$ the equation

$$\mu_t^B = \mu_0^B + \int_0^t L^B(\mu_s^B)\, ds \tag{24}$$

has a unique solution which is a continuous map $\mu_{(\cdot)}^B : t \to \mathcal{M}_B$. In fact, writing

$$\mu_t^B = \sum_{i\leq B, j\in\{0,1\}} c_{t,i,j}^B \delta_{(i,j)}$$

results in a finite approximation of the two-site Smoluchowski system:

$$\dot{c}_{t,i,j}^B = \tfrac{1}{2}\sum_{k+m=i} kmc_{t,k,j}^B c_{t,m,j}^B - ic_{t,i,j}^B \sum_{k<B} kc_{t,k,j}^B - ic_{t,B,0}^B c_{t,i,j}^B$$

$$+ \kappa(c_{t,i,1-j}^B - c_{t,i,j}^B) \quad \text{for} \quad i < B, \tag{25}$$

$$\dot{c}_{t,B,0}^B = \tfrac{1}{2}\sum_{\substack{k+m\geq B\\k,m<B}} (k+m)kmc_{t,k,j}^B c_{t,m,j}^B + c_{t,B,0}^B \sum_{k<B} k^2 c_{t,k,j}^B.$$

This system of $B$ first-order ODEs obviously has a unique global solution. We may apply Theorem 5 in the same way as was done in the first part of the proof of Theorem 2 to see that the solution is nonnegative. Moreover, it is not difficult to check that (25) conserves the total mass

$$\nu(\mu_0) \equiv \nu_B(\mu_0^B) := \int_{\mathbb{N}\times\{0,1\}, m<B} m\mu_0^B(d(m,i)) + \mu_0^B(\{(B,0)\})$$

being independent of $B$.

PROPOSITION 6. *The solution of the approximating system* (25) *with initial condition*

$$c^B_{0,1,0} = \frac{1}{1+\lambda}, \qquad c^B_{0,1,1} = \frac{\lambda}{1+\lambda}, \qquad c^B_{0,m,j} = 0, \qquad m > 1$$

*fulfills*

$$\sum_{i=1}^{B-1} i |c^\lambda_{t,i,j} - c^B_{t,i,j}| \leq c^B_{t,B,0} \leq B^{-1} C_0, \qquad 0 \leq t \leq T < t^\lambda_{\text{gel}},$$

*where $C_0$ is a constant depending on $T$ and $\lambda$.*

PROOF. Let us recall the two-site modified Smoluchowski equation (omitting the upper index $\lambda$ for simplicity)

(26)
$$\dot{c}_{t,i,j} = \tfrac{1}{2} \sum_{k+m=i} km c_{t,k,j} c_{t,m,j} - i \zeta_{t,*,i} c_{t,i,j} - \kappa c_{t,i,j} + \kappa c_{t,i,1-j},$$
$$\dot{\zeta}_{t,*,j} = -\kappa \zeta_{t,*,j} + \kappa \zeta_{t,*,1-j},$$

which we considered for the initial condition

$$c_{0,1,0} = \zeta_{0,*,0} = \frac{1}{1+\lambda}, \qquad c_{0,1,1} = \zeta_{0,*,1} = \frac{\lambda}{1+\lambda}, \qquad c_{0,m,j} = 0, \qquad m > 1.$$

The solution obeys $\zeta_{t,*,j} = \sum_k k c_{t,k,j}$ up to $t_{\text{gel}}$ as was shown in the proof of Theorem 2. Moreover, using the quantity

$$\zeta^{(B)}_{t,j} = \sum_{k \geq B} k c_{t,k,j} = \zeta_{t,*,j} - \sum_{k < B} k c_{t,k,j},$$

we may infer for $t < t_{\text{gel}}$ from (26)

(27)
$$\dot{c}_{t,i,j} = \tfrac{1}{2} \sum_{k+m=i} km c_{t,k,j} c_{t,m,j} - i c_{t,i,j} \sum_{k<B} k c_{t,k,j} - i c_{t,i,j} \zeta^{(B)}_{t,j}$$
$$\qquad - \kappa c_{t,i,j} + \kappa c_{t,i,1-j} \qquad \text{for } i < B,\ t < t_{\text{gel}},$$
$$\dot{\zeta}^{(B)}_{t,j} = \tfrac{1}{2} \sum_{\substack{k+m \geq B \\ k,m < B}} (k+m) km c_{t,k,j} c_{t,m,j} + \zeta^{(B)}_{t,j} \sum_{k<B} k^2 c_{t,k,j}$$
$$\qquad - \kappa \zeta^{(B)}_{t,j} + \kappa \zeta^{(B)}_{t,1-j}.$$

Define

$$\beta^B_{t,i,j} := c_{t,i,j} - c^B_{t,i,j}, \qquad i < B.$$

Subtracting (25) from (27) we get

$$\dot\beta^B_{t,i,j} = \tfrac{1}{2} \sum_{k+m=i} km c_{t,k,j} c_{t,m,j} - \tfrac{1}{2} \sum_{k+m=i} km c^B_{t,k,j} c^B_{t,m,j}$$

$$- i c_{t,i,j} \sum_{k<B} k c_{t,k,j} + i c^B_{t,i,j} \sum_{k<B} k c^B_{t,k,j}$$

$$- i c_{t,i,j} \zeta^{(B)}_{t,j} + i c^B_{t,i,j} c^B_{t,B,0} + \kappa(\beta^B_{t,i,1-j} - \beta^B_{t,i,j})$$

$$= \tfrac{1}{2} \sum_{k+m=i} km(\beta^B_{t,k,j} \beta^B_{t,m,j} + 2 c^B_{t,k,j} \beta^B_{t,m,j})$$

$$- i \sum_{k<B} k(\beta^B_{t,i,j} c_{t,k,j} + c^B_{t,i,j} \beta^B_{t,k,j})$$

$$- i(\beta^B_{t,i,j} \zeta^{(B)}_{t,j} + c^B_{t,i,j} \zeta^{(B)}_{t,j} - c^B_{t,B,0} c^B_{t,i,j}) + \kappa(\beta^B_{t,i,1-j} - \beta^B_{t,i,j}).$$

The mass conservation law for (25) and (27) implies

$$c^B_{t,B,0} = 1 - \sum_{l=0}^{1} \sum_{k<B} k c^B_{t,k,l} = \zeta^{(B)}_{t,j} + \zeta^{(B)}_{t,1-j} + \sum_{l=0}^{1} \sum_{k<B} k(c_{t,k,l} - c^B_{t,k,l}).$$

Thus we obtain

$$\dot\beta^B_{t,i,j} = \tfrac{1}{2} \sum_{k+m=i} km(\beta^B_{t,k,j} \beta^B_{t,m,j} + 2 c^B_{t,k,j} \beta^B_{t,m,j}) + \kappa(\beta^B_{t,i,1-j} - \beta^B_{t,i,j})$$

$$- i \sum_{k<B} k(\beta^B_{t,i,j} c_{t,k,j} + c^B_{t,i,j} \beta^B_{t,k,j})$$

$$- i\left(\beta^B_{t,i,j} \zeta^{(B)}_{t,j} - \left(\zeta^{(B)}_{t,1-j} + \sum_{k<B} k(\beta^B_{t,k,j} + \beta^B_{t,k,1-j})\right) c^B_{t,i,j}\right)$$

$$= \tfrac{1}{2} \sum_{k+m=i} km(\beta^B_{t,k,j} \beta^B_{t,m,j} + 2 c^B_{t,k,j} \beta^B_{t,m,j}) + \kappa(\beta^B_{t,i,1-j} - \beta^B_{t,i,j})$$

$$- i \sum_{k<B} k \beta^B_{t,i,j} c_{t,k,j}$$

$$- i\left(\beta^B_{t,i,j} \zeta^{(B)}_{t,j} - \left(\zeta^{(B)}_{t,1-j} + \sum_{k<B} k \beta^B_{t,k,1-j}\right) c^B_{t,i,j}\right).$$

Finally we arrive at the following finite system of ODEs:

$$\dot\beta^B_{t,i,j} = \tfrac{1}{2} \sum_{k+m=i} km(\beta^B_{t,k,j} \beta^B_{t,m,j} + 2 c^B_{t,k,j} \beta^B_{t,m,j}) + \kappa(\beta^B_{t,i,1-j} - \beta^B_{t,i,j})$$

$$- i \beta^B_{t,i,j}\left(\zeta^{(B)}_{t,j} + \sum_{k<B} k c_{t,k,j}\right) + i c^B_{t,i,j}\left(\zeta^{(B)}_{t,1-j} + \sum_{k<B} k \beta^B_{t,k,1-j}\right),$$

with initial condition $\beta^B_{0,i,j} = 0$, which has a unique solution. Again we may apply Theorem 5 in the same way as was done in the first part of the proof of Theorem 2 to see that the solution is nonnegative, taking into account the boundedness of the total mass of the solution to (27) by 1 and the fact that this solution was already proved to be nonnegative. Hence we have the relation

(28) $$c^B_{t,i,j} \le c_{t,i,j}, \qquad i < B, \ t \ge 0.$$

Thus
$$\sum_{i=1}^{B-1} |c_{t,i,j} - c^B_{t,i,j}| = \sum_{i=1}^{B-1} (c_{t,i,j} - c^B_{t,i,j})$$

and since both systems conserve the total mass we get
$$\sum_{i=1}^{B-1} i|c_{t,i,j} - c^B_{t,i,j}| = c^B_{t,B,0} - \zeta^{(B)}_{t,0} - \zeta^{(B)}_{t,1} \le c^B_{t,B,0}.$$

So, in order to estimate the total deviation of the solution of (25) with respect to that of (26) for masses less than $B$ it is enough to estimate $c^B_{t,B,0}$. We get from (25) and (28) for any $t \le T < t_{\text{gel}}$

$$\begin{aligned}
\dot c^B_{t,B,0} &= \tfrac{1}{2} \sum_{\substack{k+m \ge B \\ k,m < B}} (k+m) km c^B_{t,k,j} c^B_{t,m,j} + c^B_{t,B,0} \sum_{k < B} k^2 c^B_{t,k,j} \\
&\le \tfrac{1}{2} \sum_{\substack{k+m \ge B \\ k,m < B}} (k+m) km c_{t,k,j} c_{t,m,j} + c^B_{t,B,0} \sum_{k < B} k^2 c_{t,k,j} \\
&\le \tfrac{1}{2} \sum_{\substack{k+m \ge B \\ k,m < B}} (k+m) km c_{t,k,j} c_{t,m,j} + c^B_{t,B,0}(\sigma_{t,0} + \sigma_{t,1}) \\
&\le \tfrac{1}{2} \sum_{\substack{k+m \ge B \\ k,m < B}} (k+m) km c_{t,k,j} c_{t,m,j} + c^B_{t,B,0} C,
\end{aligned}$$

where $C$ is a constant depending on $T$ as it was shown in the proof of Theorem 2 that the second moments $\sigma_{t,j} = \sigma^\lambda_{t,j}$ are continuous functions in each compact subinterval of $[0, t_{\text{gel}})$. Hence in $[0,T]$ we have (using the notation $[x]$ for the integer part of a real number $x$, and $\zeta_{t,*,0} + \zeta_{t,*,1} \equiv 1$ for the fourth line)

$$\dot c^B_{t,B,0} \le \sum_{\substack{k+m \ge B \\ k,m < B}} k^2 m c_{t,k,j} c_{t,m,j} + c^B_{t,B,0} C$$

$$\leq \sum_{k+m\geq B} k^2 m c_{t,k,j} c_{t,m,j} + c_{t,B,0}^B C$$

$$\leq \sum_{\substack{k\geq [B/2]\\m}} k^2 m c_{t,k,j} c_{t,m,j} + \sum_{\substack{m\geq [B/2]\\k}} k^2 m c_{t,k,j} c_{t,m,j} + c_{t,B,0}^B C$$

$$\leq \sum_{k\geq [B/2]} k^2 c_{t,k,j} + (\sigma_{t,0} + \sigma_{t,1}) \sum_{m\geq [B/2]} m c_{t,m,j} + c_{t,B,0}^B C$$

$$\leq \sum_{k\geq [B/2]} k^2 c_{t,k,j} + C \sum_{m\geq [B/2]} m c_{t,m,j} + c_{t,B,0}^B C$$

$$\leq [B/2]^{-1} \left( \sum_{k\geq [B/2]} k^3 c_{t,k,j} + C \sum_{m\geq [B/2]} m^2 c_{t,m,j} \right) + c_{t,B,0}^B C$$

$$\leq [B/2]^{-1} \left( \sum_k k^3 c_{t,k,j} + C(\sigma_{t,0} + \sigma_{t,1}) \right) + c_{t,B,0}^B C$$

$$\leq [B/2]^{-1} \left( \sum_k k^3 c_{t,k,j} + C^2 \right) + c_{t,B,0}^B C$$

$$\leq [B/2]^{-1} (C' + C^2) + c_{t,B,0}^B C$$

as, again in the proof of Theorem 2, it was shown that also the third moment is bounded uniformly in compact subintervals of $[0, t_{\text{gel}})$. Hence we see that the positive quantity $c_{t,B,0}^B$ in $[0,T]$ is bounded (in view of Theorem 5) by the solution of the differential equation

$$\dot{x} = [B/2]^{-1}(C' + C^2) + Cx$$

with initial condition $x(0) = 0$, which is given by

$$x(t) = [B/2]^{-1} C^{-1} (C' + C^2)(e^{Ct} - 1).$$

Therefore, there is a constant $C_0$ depending on $T$ such that $c_{t,B,0}^B \leq B^{-1} C_0$ in $[0,T]$. This proves the proposition. $\square$

Consider the process (22) and the corresponding coupled family of approximations constructed above. Set

$$\widetilde{X}_t = N^{-1} X_{t/N}, \qquad \widetilde{X}_t^B = N^{-1} X_{t/N}^B.$$

Let $d$ be a metric on the set of all finite measures on $\mathbb{N} \times \{0,1\}$ which generates the weak topology, that is, $\langle f, \mu_n \rangle \to \langle f, \mu \rangle$ if and only if $d(\mu_n, \mu) \to 0$ for each bounded function on $\mathbb{N} \times \{0,1\}$. We may choose $d$ in a way that $d(\mu, \mu') \leq \|\mu - \mu'\|$ for each pair of measures, where $\|\cdot\|$ denotes the total variation norm.

PROPOSITION 7. *Assume*

(29) $$N_0/N \to \frac{1}{1+\lambda}, \qquad N_1/N \to \frac{\lambda}{1+\lambda} \qquad as\ N \to \infty.$$

*Then, for all* $t \geq 0$,

$$\sup_{s \leq t} d(\widetilde{X}_s^B, \mu_s^B) \to 0$$

*in probability as* $N \to \infty$.

PROOF. The sequence of laws corresponding to $(\widetilde{X}_t^B)_{t \geq 0}$ is tight, and any weak limit point (with respect to $N$) is concentrated on the set of solutions of the equation (24) (cf., e.g., [4], Theorem 2.3). Since the solution $(\mu_t^B)$ of this equation is unique, we get convergence in distribution, which implies the statement of the proposition (cf., e.g., [18], Section 3). □

The following theorem and its proof adapt Theorem 4.4 of [14] to our situation. It is a stronger version of Theorem 1.

THEOREM 8. *Assume* (29) *and consider the solution* $(\mu_t)_{0 \leq t < t_{\text{gel}}^\lambda}$ *of the Smoluchowski equation up to* $t_{\text{gel}}^\lambda$. *Then for all* $t < t_{\text{gel}}^\lambda$

$$\sup_{s \leq t} d(x\widetilde{X}_s(d[x,j]), x\mu_s(d[x,j])) \to 0$$

*in probability, as* $N \to \infty$.

PROOF. For simplicity, we omit the upper index $\lambda$ and write $dx$ instead of $d[x,j]$. Fix $\delta$ and $t < t_{\text{gel}}$. By Proposition 6 we find some $B$ such that $\mu_t^B(\{(B,0)\}) < \delta/4$. Since the measures are supported on a compact set, Proposition 7 implies

$$\sup_{s \leq t} d(x\widetilde{X}_s^B(dx \cap [1,B)), x\mu_s^B(dx \cap [1,B))) \to 0$$

in probability as $N \to \infty$. Using that the derivative of $\mu_s^B(\{(B,0)\})$ is nonnegative, we have by Proposition 6 for $s \leq t$

$$\|x\mu_s(dx) - x\mu_s^B(dx \cap [1,B))\| \leq 2\mu_s^B(\{(B,0)\}) \leq 2\mu_t^B(\{(B,0)\}) < \delta/2.$$

On the other hand,

$$\|x\widetilde{X}_s(dx) - x\widetilde{X}_s^B(dx \cap [1,B))\| \leq 2\widetilde{X}_s^B(\{(B,0)\})$$

by the fact that the evolutions of particles smaller than $B$ coincide for $\{X_t\}$ and $\{X_t^B\}$ except for those which are hidden and aggregated at $(B,0)$ in $\{X_t^B\}$, and by conservation of the total mass. We may proceed as follows, bearing in mind that the aggregated mass is increasing in time:

$$\begin{aligned}
\|x\widetilde{X}_s(dx) &- x\widetilde{X}_s^B(dx \cap [1,B))\| \\
&\leq 2\widetilde{X}_t^B(\{(B,0)\}) \\
&\leq 2\mu_t^B(\{(B,0)\}) + 2|\widetilde{X}_t^B(\{(B,0)\}) - \mu_t^B(\{(B,0)\})| \\
&\leq \delta/2 + 2\|\widetilde{X}_t^B - \mu_t^B\|.
\end{aligned}$$

So, using the property of $d$ to be dominated by the variation distance, we obtain

$$\begin{aligned}
d(x\widetilde{X}_s&(dx), x\mu_s(dx)) \\
&\leq \|x\widetilde{X}_s(dx) - x\widetilde{X}_s^B(dx \cap [1,B))\| \\
&\quad + d(x\widetilde{X}_s^B(dx \cap [1,B)), x\mu_s^B(dx \cap [1,B))) \\
&\quad + \|x\mu_s(dx) - x\mu_s^B(dx \cap [1,B))\| \\
&\leq \delta + 2\|\widetilde{X}_t^B - \mu_t^B\| + d(x\widetilde{X}_s^B(dx \cap [1,B)), x\mu_s^B(dx \cap [1,B)))
\end{aligned}$$

which gives

$$P\bigg(\sup_{s \leq t} d(x\widetilde{X}_s(dx), x\mu_s(dx)) > \delta\bigg) \to 0$$

by Proposition 7, bearing in mind that all metrics considered here are equivalent on the finite-dimensional vector space corresponding to measures on $\{1, 2, \ldots, B\} \times \{0, 1\}$. □

**Acknowledgment.** The authors would like to thank Hans Babovsky (Ilmenau) for interesting discussions.

Technische Universität Berlin
Institut für Mathematik
Strasse des 17 Juni 136
D-10623 Berlin
Germany
e-mail: siegmund@math.tu-berlin.de

Weierstrass Institute for
Applied Analysis and Stochastics
Mohrenstrasse 39
D-10117 Berlin
Germany
e-mail: wagner@wias-berlin.de